\documentclass[times,final]{article}
    \makeatletter
    \def\ps@pprintTitle{%
    \let\@oddhead\@empty
    \let\@evenhead\@empty
    \def\@oddfoot{}%
    \let\@evenfoot\@oddfoot}
    \makeatother
    \usepackage{lineno}
    \usepackage{mathrsfs}
    \usepackage[breaklinks=true,pdftex]{hyperref}
    \modulolinenumbers[1]
    \usepackage{amsmath,bm}
    \usepackage{svg}
    \usepackage{comment}
    \usepackage{multirow}
    \usepackage[a4paper,
            bindingoffset=0.1in,
            left=1.5cm,
            right=1.5cm,
            top=1.5cm,
            bottom=1.5cm]{geometry}
    \usepackage{ulem}
    \usepackage{booktabs}
    \usepackage{dirtytalk}
    \usepackage{lineno}
    \usepackage{caption}
    \usepackage{amssymb}
    \usepackage{amsmath}
    \usepackage{savesym}
        \usepackage{algorithmic}
    \usepackage{algorithm}
    \savesymbol{AND}
    \usepackage{adjustbox}

    \usepackage{rotating}
    \usepackage{graphicx}
    \usepackage{wrapfig}
    \usepackage{amssymb}
    \usepackage{soul,color}
    \usepackage{subcaption}
    \usepackage{lineno}
\usepackage{xcolor}
\usepackage{url}
    \soulregister\cite7
    \soulregister\ref7
    \soulregister\pageref7

    \usepackage{todonotes}

    \DeclareRobustCommand{\uvec}[1]{{%
    \ifcsname uvec#1\endcsname
    \csname uvec#1\endcsname
    \else
    \bm{\hat{\mathbf{#1}}}%
    \fi}}
\mathchardef\breakingcomma\mathcode`\,
{\catcode`,=\active
  \gdef,{\breakingcomma\discretionary{}{}{}}
}

\begin{document}

\begin{center}
{\LARGE \textbf{Accelerating Simulation-Driven Optimisation of Marine Propellers Using Shape-Supervised Dimension Reduction}}\\
\vspace{0.5cm}
{\small Shahroz Khan$^{1,*}$\let\thefootnote\relax\footnote{$^*$Corresponding author. E-mail address: shahroz.khan@strath.ac.uk (S. Khan)}},
{\small Stefano Gaggero$^2$},
{\small Panagiotis Kaklis$^1$},
{\small Giuliano Vernengo$^2$},
{\small Diego Villa$^2$}
\\\vspace{0.2cm}
{\small $^1$Dept. of Naval Architecture, Ocean and Marine Engineering, University of Strathclyde, Glasgow, United Kingdom}\\
{\small $^2$Dept. of Naval Architecture, Electrical, Electronic and Telecommunication Engineering, University of Genoa, Genoa, Italy}
\end{center}

\section*{\centering Abstract}
Simulation-driven shape optimisation (SDSO) of marine propellers is often obstructed by high-dimensional design spaces stemming from its complex geometry and baseline parameterisation, which leads to the notorious curse of dimensionality. In this study, we propose using the shape-supervised dimension reduction (SSDR) approach to expedite the SDSO of marine propellers by extracting latent features for a lower-dimensional subspace. SSDR is different from other dimension reduction approaches as it utilises a shape-signature vector function, which consists of a shape modification function and geometric moments, maximising the retained geometric and physical information in the subspace. The resulting shape-supervised subspace from SSDR enables us to efficiently and effectively find an optimal design in appropriate areas of the design space. The feasibility of the proposed method is tested for the E779A propeller parameterised with 40 design parameters with the objective to maximise efficiency while reducing suction side cavitation. The results demonstrate that the shape-supervised subspace achieved an 87.5\% reduction in the original design space's dimensionality, resulting in faster optimisation convergence.
\vspace{0.2cm}\\
\textit{Keywords:} Propeller Design, Computer-Aided Design, Parametric Design, Geometric Moments, Shape Optimisation.

\section{INTRODUCTION}

Over the past decade, simulation-driven shape optimisation (SDSO) has become central to the design of marine propellers. It is widely favoured to automate design exploration for maximising the efficiency of the propeller while minimising unwanted effects such as noise, vibration, and cavitation \cite{gaggero2020reduced}. The success of SDSO in finding a global optimum largely depends on the quality of the design spaces, which ensures the production of diverse and valid designs during optimisation, leading to innovative and physically plausible solutions \cite{khan2017customer}. To enhance design diversity, designers often parameterise as many features as possible. However, this approach can make it challenging to ensure design validity and, more importantly, results in the notorious curse of dimensionality, necessitating numerous computationally expensive simulation runs.

Design space dimension reduction (DSDR) techniques, also referred to as feature extraction or embedding, manifold learning, etc., depending on the field, are commonly employed to mitigate this curse. Widely used DSDR approaches include Karhunen–Loève Decomposition \cite{R07} (closely related to Principal Component Analysis (PCA) and proper orthogonal decomposition \cite{gaggero2020reduced}) and their nonlinear extensions, such as kernel PCA \cite{d2017nonlinear} and autoencoders \cite{masood2021machine}. These methods aim to extract latent variables from the original design space, creating a subspace with significantly lower dimensionality while retaining similar geometric variance. As the subspace has fewer dimensions, it requires fewer simulation runs, thus expediting the SDSO process.

However, traditional DSDR approaches often fail to preserve the complexity and intrinsic underlying geometric structure of a design. Consequently, the resulting subspace may not produce diverse and valid shapes during shape optimisation, hampering the success of SDSO by consuming the majority of the available computational budget on exploring infeasible, practically invalid, and similar designs \cite{khan2022shape}. It is important to note that these drawbacks are not necessarily intrinsic to the methods but primarily result from the low-level shape discretisations used in subspace learning.

Khan et al. \cite{khan2022shape} recently proposed a shape-supervised dimension reduction (SSDR) approach that leverages the compact geometric representation of geometries in the original design space using a shape-signature vector (SSV) function. The SSV, constructed with the shape modification function and the geometric moments of the shape's integral properties (i.e., geometric moments and their invariants), aims to maximise the accumulation of geometric variation \cite{KHAN2023116051,khan2022geometric,khan2023shiphullgan}. The SSDR then evaluates the Karhunen–Loève Expansion (KLE) of the SSV, where the solution of a variational problem allows the evaluation of latent features as a linear combination of original designs. The features provided by the KLE are expressed by the eigenfunctions of a symmetric and positive definite covariance function constructed with the SSV. The KL values associated with each feature enable the identification of active and inactive features, with active features parameterising the shape and serving as a new basis to span the subspace, maintaining the highest variance in geometry.

Khan et al. \cite{khan2022shape} validated their approach for SDSO of ship hulls, demonstrating that the resulting subspace not only had an enhanced representation capacity and compactness for producing valid and diverse design alternatives but was also physically informed to improve the convergence rate of the shape optimiser toward an optimal solution due to the correlation between hull physics (i.e., wave-making resistance) and the geometric moments present in the SSV.

The aim of this work is to explore further applications of SSDR in the SDSO of marine propellers. Additionally, this study seeks to analyse whether SSDR can still create an efficient subspace for providing diverse and valid designs while expediting the convergence of optimisation, even in the absence of known correlations between propeller physics and geometric moments. Moreover, we also compare the results of SSDR with the baseline KLE/PCA both in terms of the quality of the subspaces to provide valid and optimised geometries.  

To achieve these objectives, we first develop a multi-objective SDSO framework to optimise the shape of a marine propeller, specifically the INSEAN E779 model \cite{salvatore2006description}, with the aim of increasing propulsive efficiency and reducing cavitation occurrences. We then construct a 40-dimensional design space ($\mathcal{T}$) using the primary geometrical characteristics of the propeller, including pitch, chord, maximum and sectional camber distributions along the blade radius. Each distribution is represented by a B-spline defined through a set of control points, which represent specific dimensions of the design space. The locations of these control points determine the dimensionality and extent of the input space.

We then implemented SSDR and KLE on this design space, where for SSDR the geometric moments of the third order are used for the construction of SSV. 
Implementing SSDR and KLE resulted in 6- and 5-dimensional (6D and 5D) subspaces, achieving an 88\% and 86\% reduction in the dimensionality of the original design while preserving 95\% of the geometric variance. After constructing the 5D and 6D subspaces, they are connected to the optimiser and hydrodynamic solver. The hydrodynamic performance of the propellers is computed using a previously validated BEM (see \cite{bertetta2012cpp,gaggero2016design}). A Genetic Algorithm is implemented in the modeFRONTIER\footnote{https://www.esteco.com} optimisation environment to guide design exploration towards global optima within the subspace.

\section{SELECTED TEST CASE: E779A PROPELLER}

The selected reference test case is the INSEAN-E779A, a four-bladed, fixed-pitch, right-handed propeller, as shown in Figure 1. Originally designed in 1959, the propeller features minimal blade skew and rake, and the pitch ratio remains nearly constant along the radius. Detailed geometrical data can be found in \cite{salvatore2006description}, and the main characteristics are summarised in Table \ref{table_1}.

\begin{table}[htb!]
    \small
    \centering
    \caption{Principal characteristics of the INSEAN E779A propeller and the design point for the optimisation process.}
    \begin{tabular}{lclc}
    \toprule
    Characteristics & Value & Characteristics & Value\\
    \midrule
    Model Diameter (m) & 0.2272 & $r_{\text{hub}}/R$ & 0.22\\
    Number of blades & 4 & $A_E/A_o$ & 0.647\\
    $c/D_{\text{root}}$ & 0.256 & $c/D_{o.7R}$ & 0.378\\
    $P/D_{\text{root}}$ & 1.11 & $P/D_{0.7R}$ & 1.11\\
    Skew at tip (deg.) & 4.9 & rps (model scale) &36\\
    $J_{\text{design}}$ & 0.83 & $\sigma_{\text{design}}$ & 1.38\\
    \bottomrule
    \end{tabular}
    \label{table_1}
\end{table}

Although the INSEAN-E779A propeller has an obsolete geometry, typical of early 1960s designs, it still poses a challenging test case for validating numerical codes and associated cavitation models. Thanks to the extensive experimental datasets provided in \cite{pereira2002experimental}, this propeller has been used as a benchmark in numerous recent studies \cite{salvatore2009propeller,bensow2010simulating,vaz2015cavitating}.

The relatively simple design of the INSEAN-E779A propeller, along with the occurrence of cavitation phenomena—primarily characterised by sheet cavitation from the leading edge—makes it well-suited for the assumptions and capabilities of the Boundary Element Method (BEM) employed in this study. Furthermore, it facilitates the description of SSDR approaches, which is the primary objective of this paper.

Additionally, since there is no specific propulsive problem to address, geometric constraints such as thickness or pitch unloading at the tip can be disregarded. This expands the design space and emphasises the flexibility of the SSDR representation for reduced design space reconstruction.

    \begin{figure}[htb!] 
    \centering
    \includegraphics[width=0.8\textwidth]{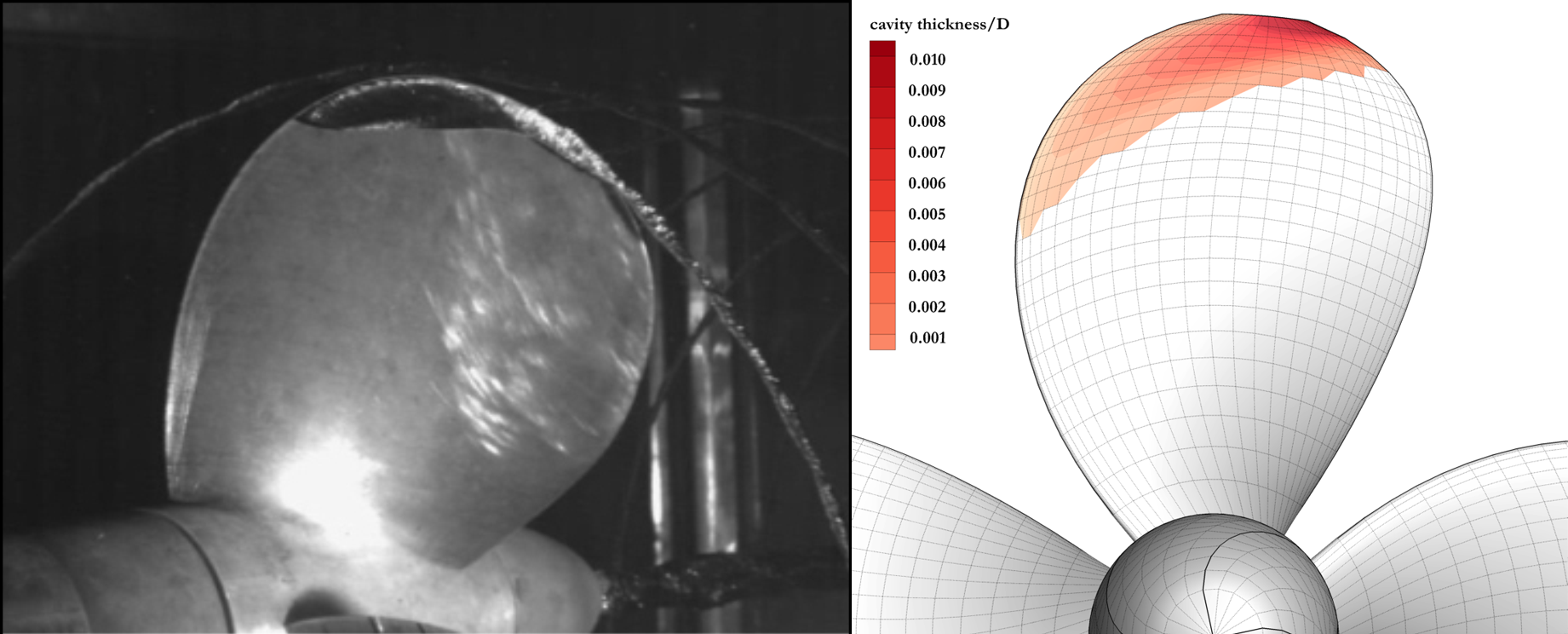}
    \caption{Predicted cavitation extension on the suction side of the E779A propeller. Comparison with experiments available in \cite{vaz2015cavitating} at the design point ($J = 0.833$, $\sigma_N = 1.38$).}
    \label{E779A_Ref}
    \end{figure}
    
An example of calculation using the Boundary Element Method is given in Fig.\ref{E779A_Ref}. By using the same $51 \times 26$ discretisation of the propeller blade adopted for the subspace construction (corresponding to $50$ chordwise panels for any of the $25$ radial strips), we computed the propeller performances at the functioning point selected for the design (advance coefficient $J = 0.833$, cavitation index based on propeller revolutions $\sigma_N = 1.38$) which represent the baseline reference for the optimisation activity. Compared to experiments \cite{salvatore2009propeller}, both delivered thrust and propulsive efficiency are well predicted. The calculated thrust coefficient is equal to $0.1761$, which is only $1.8\%$ underestimated with respect to measurements. Predicted efficiency is $2\%$ lower than experiments as the combination of the underestimation of thrust and a slight overprediction of torque. In terms of cavitation, BEM calculations predict a sheet cavity bubble from the blade leading edge that is overestimated both in radial (starting at $r/R = 0.68$ compared to $r/R = 0.85$ of experiments) and chordwise extension. These results, in any case, are consistent with similar calculations available in the literature, also using high-fidelity flow solvers like RANSE \cite{gaggero2020reduced, vaz2015cavitating}.

\section{APPROACH TO DESIGN SPACE DIMENSION REDUCTION (DSDR)}
Let $\Gamma$ be a 3D body bounded by a closed 2D manifold $\mathcal{G}$, representing a baseline/parent design of the propeller, and $\bar{\bm{\vartheta}}\in\mathcal{G}\subseteq\mathbb{R}^{\bar{n}}$, with ${\bar{n}}=1,2,3$, a coordinate set on this manifold. For an automatic shape modification, $\mathcal{G}$ is commonly parameterised with $n$ geometric parameters, defining the parametric/design vector $\mathbf{t}=(t_1,t_2,\ldots,t_n) \in \mathcal{T} \subseteq \mathbb{R}^n$. Here, $\mathcal{T}$ is the $n-$dimensional design space, which is bounded by appropriately defined set constraints, e.g., $\mathcal{T} := \left\{ \mathbf{t}: t_i^l \leq t_i \leq t_i^u, \forall i \in \lbrace 1,2,\dots n \rbrace \right\}$ with $\mathbf{t}^l,\mathbf{t}^u\in\mathbb{R}^n$ denoting the lower and upper bound vector, respectively. The parametric vector $\mathbf{t}$ of $\mathcal{G}$ yields a continuous shape modification vector $\bm{G}(\bar{\bm{\vartheta}},\mathbf{t})\in\mathbb{R}^{n_G}$ with ${n_G}=1,2,3$, which for any  $\mathbf{t}\in\mathcal{T}$ modifies the initial $\bar{\bm{\vartheta}}$ to produce new ${\bar{\bm{\vartheta}}}^\prime$ that defines the modified $\mathcal{G}^\prime$, i.e., 
    
    \begin{equation}\label{E1}
       {\bar{\bm{\vartheta}}}^\prime= \bar{\bm{\vartheta}}+\bm{G}(\bar{\bm{\vartheta}},\mathbf{t}),\,\forall \bar{\bm{\vartheta}}.
    \end{equation}

\subsection{Shape-supervised dimension reduction (SSDR)}
     To construct the subspace using SSDR, we consider that along with the continuous shape modification vectors, $\bm{G}(\bar{\bm{\vartheta}},\mathbf{t})$, there is a lumped geometric moment vector, $\bm{M}(\bm{\vartheta}_M,\mathbf{t}) \in \mathbb{R}^{n_M}$ with $n_M=1,2,\dots$, which has a null measure and corresponds to an arbitrary point, $\bm{\vartheta}_M$, where this moment vector is virtually defined. We further assume $\mathcal{G}$ and $\mathcal{M}$ as domains of definition for $\bm{G}(\bar{\bm{\vartheta}},\mathbf{t})$ and $\bm{M}(\bm{\vartheta}_M,\mathbf{t})$, respectively. Now, consider a combined geometry and moment vector $\bm{P}(\bm{\vartheta} ,\mathbf{t})\in\mathbb{R}^{n_P}$, $n_P= n_G+n_M$, defined in the domain $\mathcal{P}:={\mathcal{G}\cup\mathcal{M}}$ with $\bm{\vartheta}=(\bm{\bar{\vartheta}},\bm{\vartheta}_M)$ and 
    
    \begin{equation}\label{E28}
       \bm{P}(\bm{\vartheta} ,\mathbf{t})=\left(\bm{G}(\bm{\bar{\vartheta}},\mathbf{t}),\bm{M}(\bm{\vartheta_M},\mathbf{t})\right).
    \end{equation} 
    
    \noindent $\bm{P}(\bm{\vartheta},\mathbf{t})$ contains both the geometry and its moments and forms a unique SSV function encompassing high-level information about the baseline design. Also consider that $\bm{P}(\bm{\vartheta},\mathbf{t})$ belongs to a disjoint Hilbert space $L^2_f(\mathcal{P})$, which is defined by the generalised inner product:

    \begin{equation}\label{E29}
    \begin{aligned}
        (\mathbf{a},\mathbf{b})_f=\int_{\mathcal{P}}f(\bm{\vartheta})\mathbf{a}(\bm{\vartheta})\cdot\mathbf{b}(\bm{\vartheta})\text{d}\bm{\vartheta} \\= \int_{\mathcal{G}}f(\bm{\bar{\vartheta}})\mathbf{a}(\bm{\bar{\vartheta}})\cdot\mathbf{b}(\bm{\bar{\vartheta}})\text{d}\bm{\bar{\vartheta}} + f(\bm{\vartheta}_M)\mathbf{a}(\bm{\vartheta}_M)\cdot\mathbf{b}(\bm{\vartheta}_M),
    \end{aligned}
    \end{equation}

    \noindent with the associated norm $\Vert\mathbf{a}\Vert=(\mathbf{a},\mathbf{a})^{\frac{1}{2}}_f$, where $f(\bm{\bar{\vartheta}}),f(\bm{\vartheta}_M)\in\mathbb{R}$ are appropriate positive weight functions used to focus analysis on certain regions of $\mathcal{G}$. Considering all the possible realisations of $\mathbf{t}$ coming from $\mathcal{T}$, one can evaluate the associated mean vector of $\bm{P}$ as

    \begin{equation}\label{E30}
       \langle \bm{P} \rangle=  \int_\mathcal{T}f(\bm{\vartheta})\bm{P}(\bm{\vartheta},\mathbf{t})\rho(\mathbf{t})\text{d}\mathbf{t},
    \end{equation}
     \noindent and the associated geometrical variance as 
    \begin{equation}\label{E31}
        \sigma^2 = \langle\Vert \overline{\bm{P}}\Vert^2 \rangle = \int_{\mathcal{T}}\int_{\mathcal{P}}f(\bm{\vartheta}) \overline{\bm{P}}(\bm{\vartheta},\mathbf{t})\cdot\overline{\bm{P}}(\bm{\vartheta},\mathbf{t})\rho(\mathbf{t})\text{d}\bm{\vartheta}\text{d}\mathbf{t},
    \end{equation}

    \noindent where $\overline{\bm{P}}$ is the \textit{deviation from the mean of SSV} (i.e., $\overline{\bm{P}}=\bm{P}-\langle \bm{P} \rangle$) and $\langle \cdot \rangle$ is the ensemble average over $\mathbf{t}$. The aim of KLE is to find an optimal basis of orthonormal functions for the linear representation of $\overline{\bm{P}}(\bm{\vartheta},\mathbf{t})$: 
     
    \begin{equation}\label{E7}
        \overline{\bm{P}}(\bm{\vartheta},\mathbf{t}) \approx \sum_{i=1}^mv_i\bm{\omega}_i(\bm{\vartheta}),
    \end{equation}
    
    \noindent where 
    
    \begin{equation}\label{E8}
        v_i = \left(\overline{\bm{P}},\bm{\omega}_i\right)_f=\int_{\mathcal{P}}f(\bm{\vartheta})\overline{\bm{P}}(\bm{\vartheta},\mathbf{t})\cdot\bm{\omega}_i(\bm{\vartheta})\text{d}\bm{\vartheta}, 
    \end{equation}
    
    \noindent are the basis-function components usable as new \textit{Geometrically- and Functionally-Active Latent Variable} (GFALV) vector for shape optimisation of the propeller. The optimality condition associated with the KLE refers to the geometric variance retained by the basis functions through Eq. \eqref{E7}. Therefore, combining Eq. \eqref{E31}, \eqref{E7} and \eqref{E8} yelids: 
    
    \begin{equation}\label{E9}
    \begin{aligned}
        \sigma^2=\sum_{i=1}^\infty \sum_{j=1}^\infty \langle v_iv_j\rangle (\bm{\omega}_i(\bm{\vartheta}),\bm{\omega}_i(\bm{\vartheta}))_f \\= \sum_{j=1}^\infty\left\langle v_j^2\right\rangle = \sum_{j=1}^\infty\left\langle\left(\overline{\bm{P}},\bm{\omega}_i(\bm{\vartheta})\right)^2_f\right\rangle.
    \end{aligned}
    \end{equation}
    
    \noindent The basis retaining the maximum variance is provided by the solution of the following variational problem \cite{R07}:
    
    \begin{equation}\label{E10}
        \begin{aligned}
        \min_{\bm{\omega}\in L^2_f\left(\mathcal{P}\right)} \quad & J\left(\bm{\omega}(\bm{\vartheta})\right) = \left\langle\left(\overline{\bm{P}},\bm{\omega}(\bm{\vartheta})\right)^2_f\right\rangle \\
        \textrm{subject to} \quad & \left(\bm{\omega}(\bm{\vartheta}),\bm{\omega}(\bm{\vartheta})\right)^2_f=1,\\
        \end{aligned}
    \end{equation}
   
    \noindent which, as proven in \cite{R07}, yields
   
    \begin{equation}\label{E11}
       \mathscr{L}\bm{\omega}(\bm{\vartheta})=\int_{\mathcal{P}}f(\bm{\theta})\left\langle\overline{\bm{P}}(\bm{\vartheta},\mathbf{t})\otimes \overline{\bm{P}}(\bm{\theta},\mathbf{t})\right\rangle\bm{\omega}(\bm{\theta})\text{d}\bm{\theta}=\lambda\bm{\omega}(\bm{\vartheta}),
    \end{equation}
   
    \noindent where $\otimes$ is the outer product, $\bm{\theta},\bm{\vartheta}\in\mathcal{G}$, and $\mathscr{L}$ is the self-adjoint integral operator whose eigensolutions form the basis function for the linear representation of $\overline{\bm{P}}(\bm{\theta},\mathbf{t})$ given in Eq. \eqref{E7}. The resulting eigenvectors, or KL-modes $\left\{\bm{\omega}_i(\bm{\vartheta})\right\}_{i=1}^\infty$, are orthogonal and constitute a complete basis for $L^2_f\left(\mathcal{G\cup\mathcal{M}}\right)$. Additionally, the eigenvalues or KL-values $\left\{\lambda_i\right\}_{i=1}^\infty$ represent the variance,
    
    \begin{equation}\label{E12}
        \sigma^2=\sum_{i=1}^\infty\lambda_i,
    \end{equation}
    
    \noindent retained by the associated basis. The first $m$ eigenvectors, i.e., $\left\{\bm{\omega}_i(\bm{\vartheta})\right\}_{i=1}^m$ constitute the optimal basis for the approximation in Eq.~\eqref{E7}. Moreover, considering $\varepsilon$ as the desired level of confidence for capturing the variance, $m$ in Eq. \eqref{E7} can be selected to satisfy
    
    \begin{equation}\label{E13}
        \sum_{i=1}^m\lambda_i\geq\varepsilon\sum_{i=1}^\infty\lambda_i = \varepsilon\sigma^2
    \end{equation}
    
    \noindent with $0<\varepsilon\leq1$ and $\lambda_i\geq\lambda_{i+1}$.
    
    The numerical implementation of Eq. \eqref{E11} -- or its generalised form; see Eq. \eqref{E29} -- is performed using the approach of Diez et al. in \cite{R07}. In this approach, the integral in Eq. \eqref{E29} is computed by discretising the domain of integration, $\bar{\bm{\vartheta}}\in\mathcal{G}$, into $E$ quadrilateral mesh elements with measure equal to $\Delta\mathcal{G}_i$ and centroid at $\left\{\bar{\bm{\vartheta}}_i,i=1,2,\dots,E\right\}$. We then use the spatial discretisation $\mathbf{d}(\mathbf{t})$ and $\mathbf{W}$ of $\overline{\bm{P}}(\bm{\vartheta},\mathbf{t})$ and $\bm{\omega}(\bm{\vartheta})$, respectively. Finally, the problem is recast as an eigenproblem of a matrix $(\mathbf{A})$:
    
    \begin{equation}\label{E21}
        \mathbf{A}\mathbf{W}=\mathbf{W}\bm{\Lambda},
    \end{equation}

    \noindent where, $\mathbf{W}=\left\{ \mathbf{w}^i,i=1,2,\dots,n_GE+n_M \right\}$ is a square matrix whose $i\text{th}$ column, $\mathbf{w}^i$, is the corresponding eigenvector or KL-mode. The KL-values, $\bm{\Lambda} = \left\{\lambda_i,i=1,2,\dots,n_GE+n_M\right\}$, represent the variance retained by the associated KL-mode. For example, at $n_P=4$ (with $n_G=3$ and $n_M=1$), $\mathbf{A}$ can be represented as
   
    \begin{equation}\label{E22}
       \mathbf{A} = 
       \begin{bmatrix}
       {\mathbf{C}}_{11} & {\mathbf{C}}_{12} & {\mathbf{C}}_{13} & \mathbf{C}_{14}\\
       {\mathbf{C}}_{12} & {\mathbf{C}}_{22} & {\mathbf{C}}_{23}&
       \mathbf{C}_{24}\\
       {\mathbf{C}}_{13} & {\mathbf{C}}_{32} & {\mathbf{C}}_{33} &
       \mathbf{C}_{34}\\
       \mathbf{C}_{14} & \mathbf{C}_{24} & \mathbf{C}_{34} & \mathbf{C}_{44}   
       \end{bmatrix}
       \begin{bmatrix}
        {\mathbf{Q}} & 0 & 0 & 0\\
        0 & {\mathbf{Q}} & 0 & 0\\
        0 & 0 & {\mathbf{Q}} & 0\\
        0 & 0 & 0 & \mathbf{Q}\\
        \end{bmatrix},
    \end{equation}
    
    \noindent where $\mathbf{C}_{lk}=\left\langle\text{d}_{l}(\mathbf{t})~\left[ \text{d}_{k}(\mathbf{t})\right]^T\right\rangle$, $\forall~l,k=1,2,\dots,n_P$ and $\mathbf{Q}$ is the weighted matrix to normalise $\mathbf{C}_{lk}$, so all of its components have same influence while computing $\mathbf{A}$. For dimensionality reduction we first rearrange KL-values in $\bm{\Lambda}$ in descending order, i.e., $\lambda_i\geq\lambda_{i+1}$. Afterwards, we select the first $m$ KL-values $\left\{\lambda_i\right\}_{i=1}^m$ via Eq. \eqref{E13} along with their associated KL-modes $\left\{\mathbf{w}^i\right\}_{i=1}^m$, which correspond to features with the greatest impact on geometry changes. The spatial discretisation of $\overline{\bm{P}}(\bm{\vartheta},\mathbf{t})$ and $\bm{\omega}(\bm{\vartheta})$ (namely $\mathbf{d}(\mathbf{t})$ and $\mathbf{W}$) can now be approximated and defined as 
 
    \begin{equation}\label{E19}
        \mathbf{d}(\mathbf{t}) = 
        \begin{Bmatrix}
        \overline{P}_1(\bar{\bm{\vartheta}}_1,\mathbf{t})\\\vdots\\\overline{P}_1(\bar{\bm{\vartheta}}_E,\mathbf{t})\\\overline{P}_2(\bar{\bm{\vartheta}}_1,\mathbf{t})\\\vdots\\\overline{P}_2(\bar{\bm{\vartheta}}_E,\mathbf{t})\\\overline{P}_3(\bar{\bm{\vartheta}}_1,\mathbf{t})\\\vdots\\\overline{P}_3(\bar{\bm{\vartheta}}_E,\mathbf{t})\\\overline{P}_1(\bm{\vartheta}_M,\mathbf{t})
        \end{Bmatrix}\approx\sum_{i=1}^{m}v_i\mathbf{w}^i;~~
        \mathbf{w}^i = 
        \begin{Bmatrix}
        \omega_1(\bar{\bm{\vartheta}}_1)\\\vdots\\\omega_1(\bar{\bm{\vartheta}}_E)\\\omega_2(\bar{\bm{\vartheta}}_1)\\\vdots\\\omega_2(\bar{\bm{\vartheta}}_E)\\\omega_3(\bar{\bm{\vartheta}}_1)\\\vdots\\\omega_3(\bar{\bm{\vartheta}}_E)\\\omega_1(\bm{\vartheta}_M)
        \end{Bmatrix}.
    \end{equation}
 
    \noindent The latent variables $\bm{v}\in\mathbb{R}^m$ formulated in Eq. \eqref{E8} can be finally obtained in a discretised form as 

    \begin{equation}\label{E23}
       v_i=\mathbf{d}(\mathbf{t})^T
       \begin{bmatrix}
        \mathbf{Q} & 0 & 0 & 0\\
        0 & \mathbf{Q} & 0 & 0\\
        0 & 0 & \mathbf{Q} & 0\\
        0 & 0 & 0 & \mathbf{Q}\\
        \end{bmatrix}
        \mathbf{w}^i.
   \end{equation}
   
   \noindent It should be noted that the KL-modes are formulated while taking into account both geometry and geometric moments in order to preserve the underlying structure of $\mathcal{G}$ and to accumulate the functional information of designs in $\mathcal{T}$. Therefore, by using only the first $n_GE$ elements of column vector $\mathbf{w}^i$ in Eq. \eqref{E23} one could form the latent variable vector which is used for the shape modification of $\mathcal{G}$ during the shape optimisation performed in the subspace $\mathcal{V} := \left\{\bm{v}: v_i^l \leq v_i \leq v_i^u, \forall i \in \lbrace 1,2,\dots m \rbrace \right\}$.

The common approach for finding $v_i^l$ and $v_i^u$ employs the standard deviation from the mean shape lying at the centroid of the design space, which is evaluated as 

    \begin{equation}\label{E53}
        v_i \in \left[-\sqrt{\kappa\lambda_i},\sqrt{\kappa\lambda_i}\right],\quad\kappa\in\{1,2,3\}.
    \end{equation}

        \subsection{Geometric moments}\label{s_qm}
	In the construction of SSV, introduced in Eq. \eqref{E28}, we use a finite number of moments of  $\Gamma$, which are defined by the following equation;
	\begin{equation}\label{E35}
		M_{p,q,r} =  \int_{-\infty}^{+\infty}\int_{-\infty}^{+\infty}\int_{-\infty}^{+\infty} x^p~y^q~z^r~\rho(x,y,z)~\text{d}\Gamma,~ p,q,r\in\lbrace 0,1,2,\dots\rbrace,
	\end{equation}
    which evaluates the $s\,\text{th}$-order geometric moments of $\Gamma$, where $s=p+q+r$ and $\rho(x,y,z)=\begin{cases} 1 & \text{if}~(x,y,z)\in \Gamma\\ 0 & \text{otherwise} \end{cases}$. Given now a non-negative integer $s$, we consider the vector $\bm{M}^s$ to contain all $M_{p,q,r}$ moments for which $p + q + r = s$.  For instance, $\bm{M}^2 = \left\{  M_{2,0,0},M_{0,2,0},M_{0,0,2},M_{1,1,0},M_{1,0,1},M_{0,1,1}\right\}\in\mathbb{R}^{n_M=6}$. Furthermore, the zeroth- and first-order moments, i.e., $M_{0,0,0}$ and $M_{1,0,0}$, $M_{0,1,0}$, $M_{0,0,1}$, are commonly used in computer-aided design and engineering packages to compute an object's volume, $V=M_{0,0,0},$ and its centroid $\mathbf{c} = \left\{ C_x,C_y,C_z\right\} = \left\{ \frac{M_{1,0,0}}{M_{0,0,0}},\frac{M_{0,1,0}}{M_{0,0,0}},\frac{M_{0,0,1}}{M_{0,0,0}}\right\}$. If $\rho(x,y,z)$ is the PDF of a continuous random variable then $\bm{M}^0=1$, whereas $\bm{M}^1$, $\bm{M}^2$, $\bm{M}^3$ and $\bm{M}^4$, represent the mean, variance, skewness and kurtosis of the random variable, respectively.

     There exist several methods available in the literature for evaluating geometric moments, but the most commonly used approach is via Gauss's divergence theorem \cite{yang1997fast,M1,M2}, which allows for the conversion of volume integrals to integrals over the bounding surface(s).

    \subsection{Propeller optimisation in subspace} \label{propOpt}
    After performing the dimensionality reduction, a multi-objective optimisation, shown in Eq.\eqref{opt_eq}, of the E779A propeller is carried out. This optimisation simultaneously targets the maximisation of the open water efficiency ($\eta_o$) and the minimisation of the blade area subjected to cavitation, ($A_{\text{cavitation}}$). Specific requirements on the propeller thrust coefficient are included as optimisation constraints. More specifically the optimised propeller should be characterised by a thrust coefficient $K_T$ within $\pm1.5\%$ reference propeller value.
    
    \begin{equation}\label{opt_eq}
    \begin{aligned}
    \textrm{Find } \bm{v^*}\in\mathbb{R}^m \quad & \textrm{such that} \\
    \eta_o(\bm{v^*}) = &  \max_{\bm{v} \in \mathcal{V}} \eta_o(\bm{v}) \\
    A_{\text{cavitation}}(\bm{v^*}) = &  \min_{\bm{v} \in \mathcal{V}} A_{\text{cavitation}}(\bm{v}) \\
    \textrm{subject to} \quad  & K_{T\text{min}}\leq K_{T}(\bm{v})\leq K_{T\text{max}},
     \end{aligned}
    \end{equation}

\section{EXPERIMENTS: DIMENSION REDUCTION AND OPTIMISATION}

    This section presents the results of comprehensive experiments conducted using the proposed approach, analysing its performance for multi-objective optimisation of the E779A propeller and demonstrating its capabilities for efficient dimensionality reduction in comparison to other existing methods.
    
\subsection{Evaluation of geometric moments}
Geometric moments of any order can be calculated for geometries satisfying the conditions indicated in Section \ref{s_qm}. However, high-order geometric moments can be sensitive to noise \cite{M10} while at the same time, numerical inaccuracies are ever-present when evaluating high-order terms \cite{M7}. Furthermore, a literature review in various application areas, ranging from kinetic equations \cite{M8} to shape retrieval \cite{xu2008geometric}, reveals that moments of an order higher than four are rarely useful. We limited the order of geometric moments invariants appearing in SSV to $s=3$ in this connection, which contains $n_M=10$ components. The moment invariants for the single blade of the E779A propeller are presented in Tables \ref{table_2}. 

    \begin{table}[htb!]
    \small
    \centering
    \caption{Geometric moment invariants of $3\text{rd}-$order evaluated for the E779A propeller.}
    \begin{tabular}{lllll}
    \toprule
    $MI_{0,0,3}$ & $MI_{0,1,2}$ & $MI_{0,2,1}$ & $MI_{0,3,0}$& $MI_{1,0,2}$\\
    \midrule
     2.933E-01 &-5.375E-03&	7.557E-03&	2.131E-03	&-5.132E-03\\
     \midrule
     $MI_{1,1,1}$ & $MI_{1,2,0}$ & $MI_{2,0,1}$ & $MI_{2,1,0}$ & $MI_{3,0,0}$\\
     \midrule
     2.523E-02&	1.084E-03&	-6.153E-02	&-3.516E-03&	5.313E-03\\
    \bottomrule
    \end{tabular}
    \label{table_2}
    \end{table}

\subsection{Implementing SSDR for subspace construction}
\begin{figure}[htb!] 
    \centering
    \includegraphics[width=0.5\textwidth]{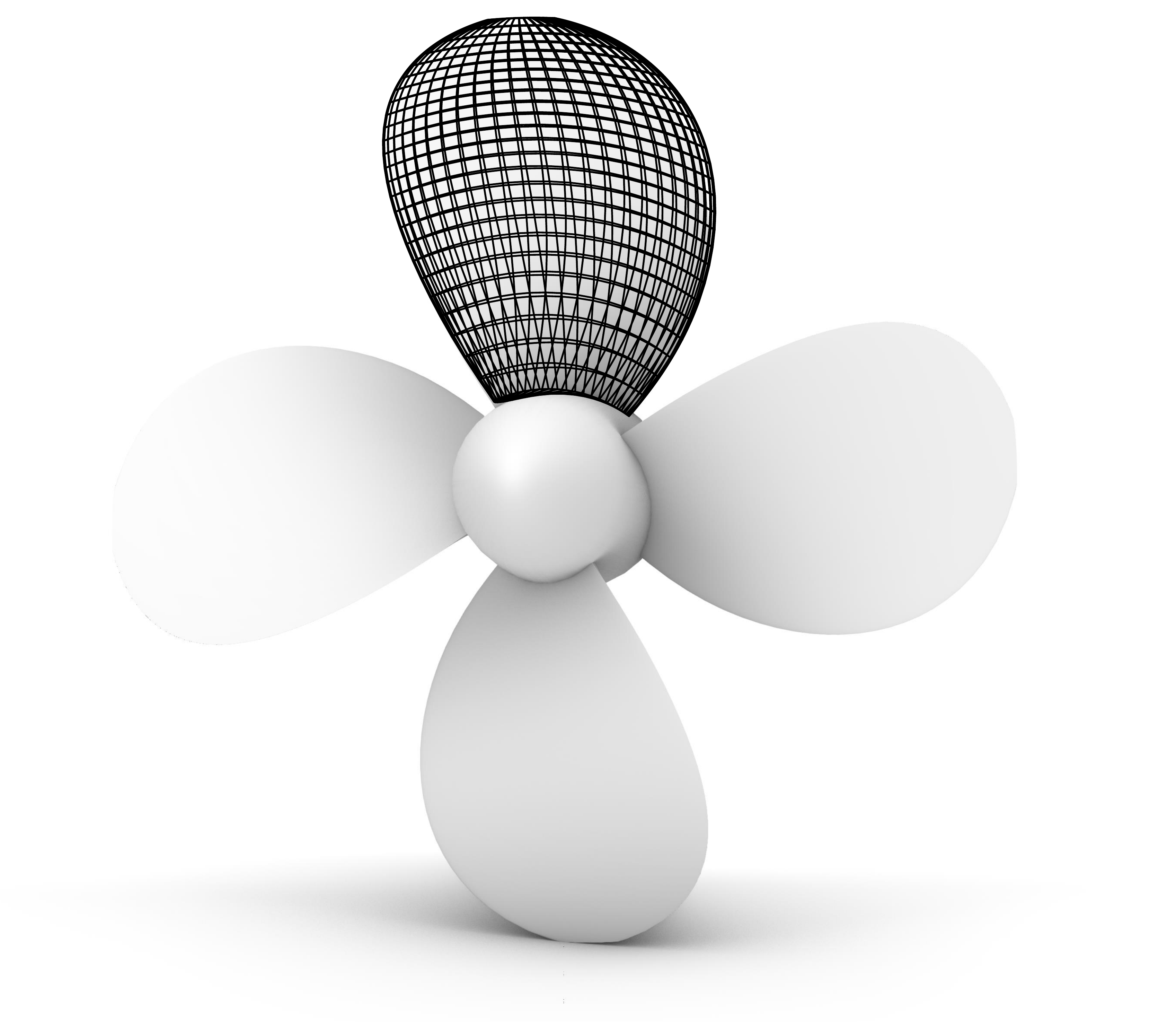}
    \caption{Computational grid of the E779A propeller used for the numerical implementation of DSDR and the BEM simulation.}
    \label{grid}
    \end{figure}

The implementation of SSDR commences with the definition of bounding limits for parameters in $\mathcal{T}$. In the present work, these bounding limits are set according to \cite{gaggero2020reduced} as they provide sufficient variation with a relatively large number of valid shapes. During SSDR, the ensemble averages, $\langle\cdot\rangle$ (in Eq. \eqref{E30}), over $\mathcal{T}$ is evaluated using Monte Carlo sampling, with statistically converged number of samples $\Psi=10000$, $\lbrace\mathbf{t}_\psi\rbrace_{\psi=1}^\Psi\sim\rho(\mathbf{t})$. $\rho(\mathbf{t})$ is selected to be a uniform distribution, thus each shape in $\mathcal{T}$ has the same possibility to be optimal. Afterwards, SSV, $\left(\bm{G}(\bm{\vartheta},\mathbf{t}),\bm{M}^3\right)$, is constructed using the third-order ($s=3$) geometric moments and the employed grid for the baseline propeller is composed of $E=51 \times 26$ nodes (see Fig. \ref{grid}), which, along with $n_G=3$ and the moments, when provided, will produce the matrix $\mathbf{A}$ in Eq. \eqref{E21} of $3988\times 3988$ dimensionality, which commence the construction of shape supervised subspace $\mathcal{V}$.

\begin{figure}[htb!] 
    \centering
    \includegraphics[width=01\textwidth]{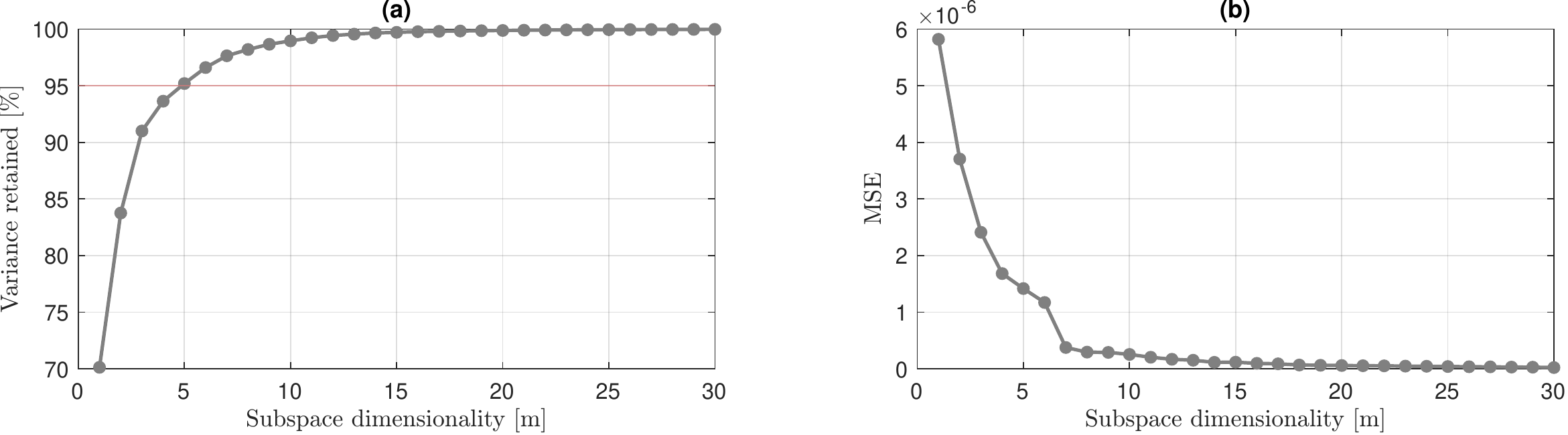}
    \caption{(a) Percentage of variance retained by each of the propeller model’s subspace versus its dimension. The horizontal red line indicates the 95\% threshold. (b) Reconstruction accuracy of propeller model’s subspaces measured via MSE with respect to their dimensionality (m).}
    \label{variance}
    \end{figure}

It should be also noted that vector spaces are normalised to exhibit the same variance associated with geometry and moment invariants. The selection of active KL-modes (eigenvectors) for the construction of subspaces is performed in a way that guarantees that every subspace retains at least $95\%$ of the variance associated with $\mathcal{T}$. In other words, the number is determined by the sum of KL-values (eigenvalues) that reach this threshold; see Eq. \eqref{E12}. Fig. \ref{variance} (a) and (b) depict the percentage of variance retained and mean square error (MSE), respectively, with respect to the dimensionality of $\mathcal{V}$ and the dimension required for each subspace to reach that level. MSE is measured between original designs in the training dataset and the reconstrued designs from the latent space of varying dimensionality. It can be seen that $\mathcal{V}$ requires $m = 5$ dimensions to capture 95\% of variance; thus, resulting in an  87.50\% of dimensional reduction, i.e., from $n = 40$ to $m = 5$. This shows the strong DSDR capability of the SSDR approach. 

\begin{figure}[htb!] 
    \centering
    \includegraphics[width=01\textwidth]{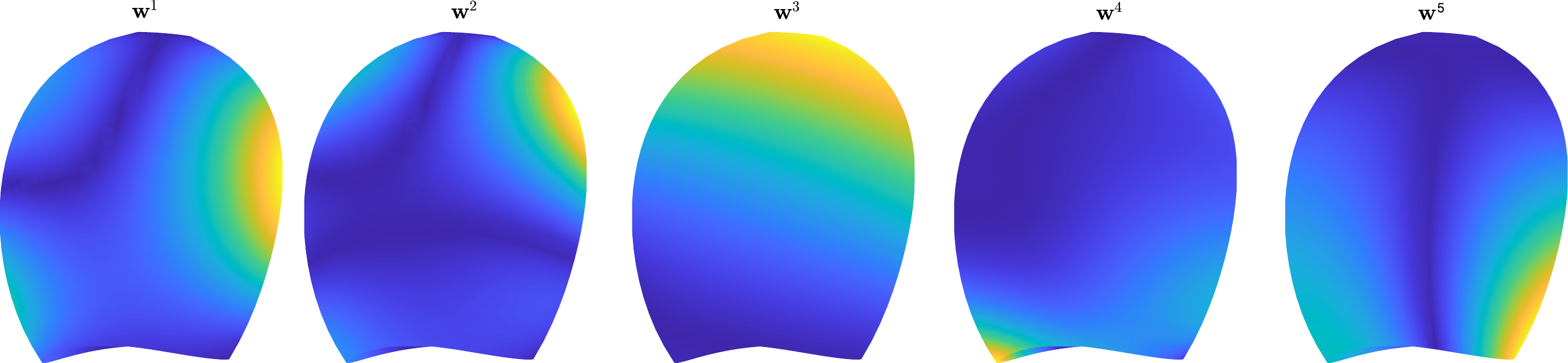}
    \caption{Shape deformation of the propeller model corresponding to the first five eigenvectors. The average magnitude of surface displacement is colour coded [small:blue to large:yellow].}
    \label{eign}
    \end{figure}

  Fig.~\ref{eign} show the first five KL-modes, $\mathbf{w}^1$, $\mathbf{w}^2$, $\mathbf{w}^3$, $\mathbf{w}^4$ and $\mathbf{w}^5$  projected on the propeller grids. This projection is of great practical value as it highlights the type and order of variance corresponding to each KL mode. From Fig.~\ref{eign} it can be seen that the first  ($\mathbf{w}^1$) and second ($\mathbf{w}^2$) KL-modes of $\mathcal{V}$ show high deviation at the trailing edge of the propeller blade. However, $\mathbf{w}^1$ this effect is more towards the root, whereas for $\mathbf{w}^2$ this effect is more towards the tip. The $\mathbf{w}^3$ show the highest deviation at the tip but towards the trailing edge and $\mathbf{w}^4$ and $\mathbf{w}^5$ exhibit deviation close to the root at the leading and trailing edges, respectively.  

\subsection{Comparison of SSDR with baseline KLE}
With the SSDR and the KLE reduced design space, we addressed the same optimisation problem. The idea, indeed, is to compare the results of the dimension reduction and the optimisation results obtained with the SSDR and with the baseline KLE in the absence of moments with the results obtained in the non-reduced design space.

\begin{figure}[htb!] 
    \centering
    \includegraphics[width=01\textwidth]{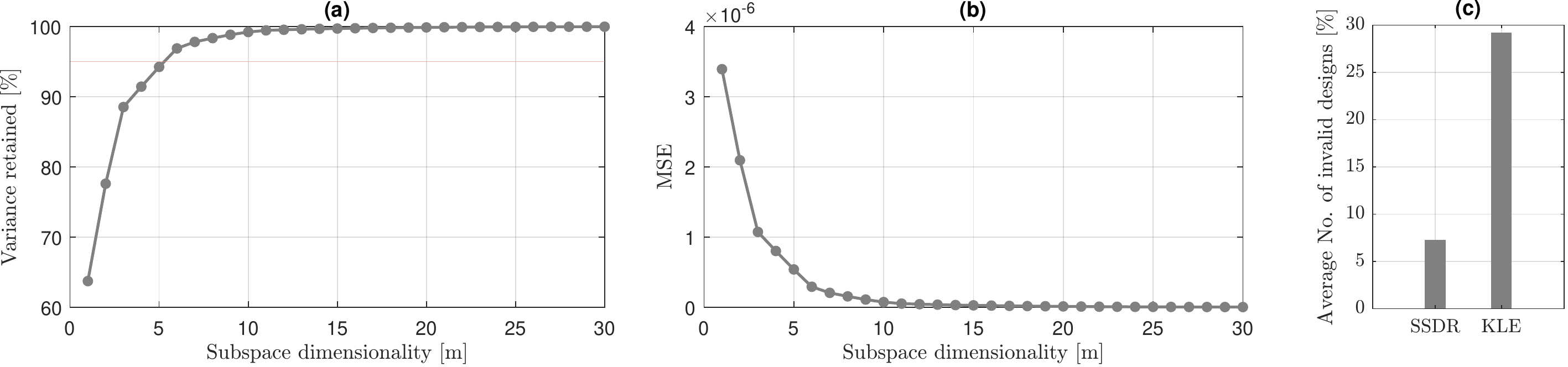}
    \caption{(a) Percentage of variance retained by each of the propeller model’s subspace obtained with KLE versus its dimension. The horizontal red line indicates the 95\% threshold. (b) Reconstruction accuracy of propeller model’s subspaces obtained with KLE measured via MSE with respect to their dimensionality (m). (c) The average per cent of invalid propeller designs obtained from the subspace generated with SSDR and KLE.}
    \label{variance_error_pca}
    \end{figure}

Similar to Fig. \ref{eign}, Fig. \ref{variance_error_pca} (a) and (b) show the percentage of variance retained and MSE versus the dimensionality of $\mathcal{V}$. It can be seen that KLE requires $m=6$ dimensions to capture 95\% of the geometric variance, which is slightly higher than the dimensionality of SSDR (i.e., $m=5$). Although in this case, the MSE is lower compared to SSDR, due to the absence of moments, the quality of KLE's subspace is significantly lower compared to SSDR in terms of generating valid designs (i.e., designs with no self-intersecting geometries). Fig. \ref{variance_error_pca} (c) shows the percentage of an average number of invalid designs resulting from the subspace created with SSDR and KLE. This quality analysis assesses the suitability of the subspace for shape optimisation, i.e., we assess whether the subspace $\mathcal{V}$ resulting from new parametrisation of shapes with latent variables $\bm{v}$ can capture the underlying shape structure adequately and whether it produces valid and diverse geometries. To commence these analyses, we use five random Monte Carlo samplings of $\Psi=5,000,000$ parameter vectors from each subspace and compute the average number of invalid shapes (i.e., shapes with self-intersecting geometries) appearing in each subspace. Fig \ref{variance_error_pca} (c) shows that the subspace resulting from SSDR produces a significantly lower number of geometries compared to the subspace of KLE. These results confirm the ability of SSDR to generate a large number of valid geometries, thereby promoting fast convergence in optimisation. This quality of the SSDR is mainly due to the inclusion of the geometric moments, which ensure the provision of a rich geometric representation of the propeller during dimension reduction by capturing all the intrinsic geometric features. 

    \subsection{Shape optimisation}
In this subsection, we discuss the results of the shape optimisation defined in the section \ref{propOpt}. This includes the optimisation performed in the original 40-dimensional non-reduced design space and the optimisation performed in the reduced dimensional subspace resulting from the baseline KLE and SSDR approaches.

    \begin{figure}[htb!] 
    \centering
    \includegraphics[width=0.48\textwidth]{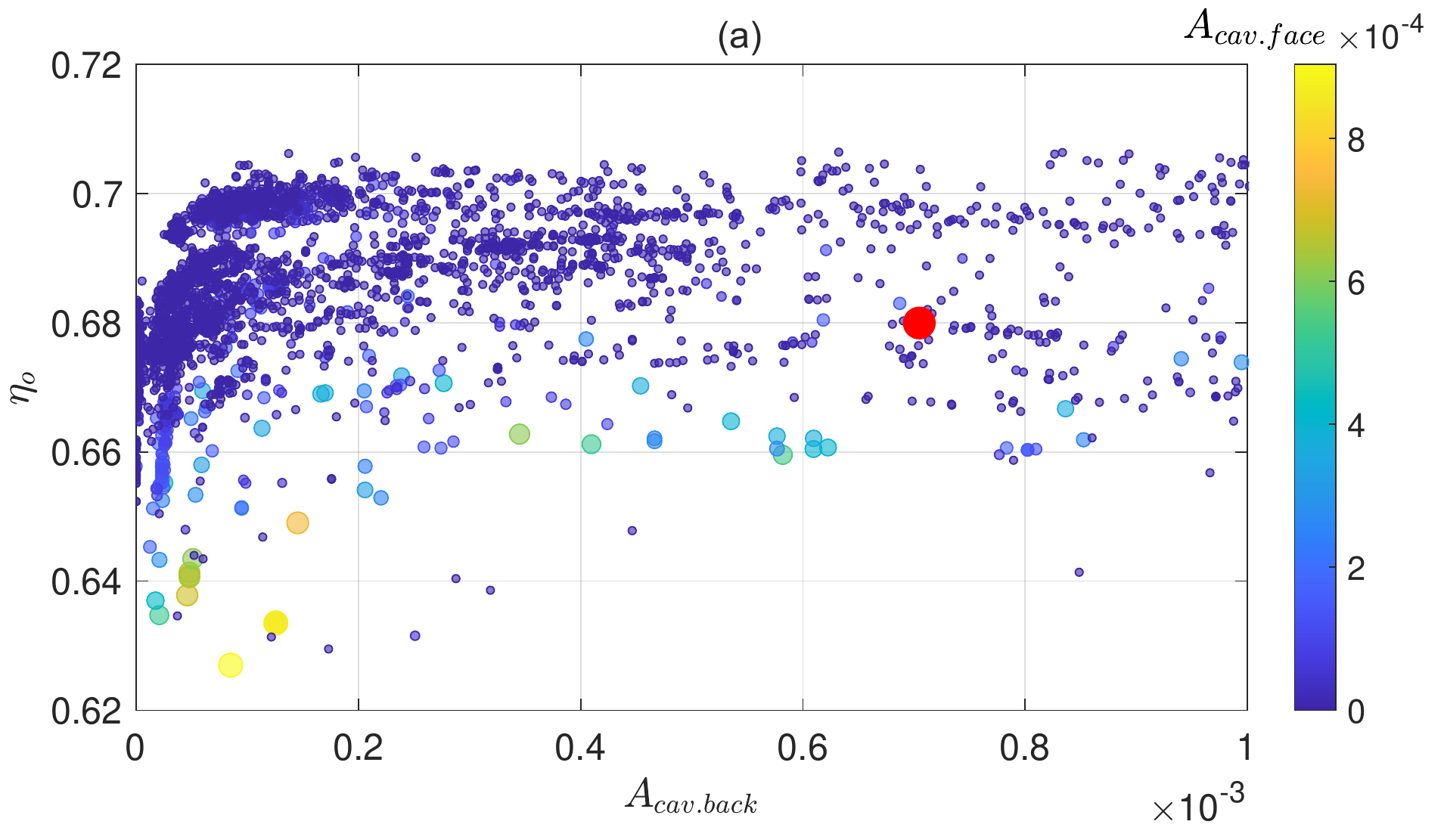}
    \includegraphics[width=0.48\textwidth]{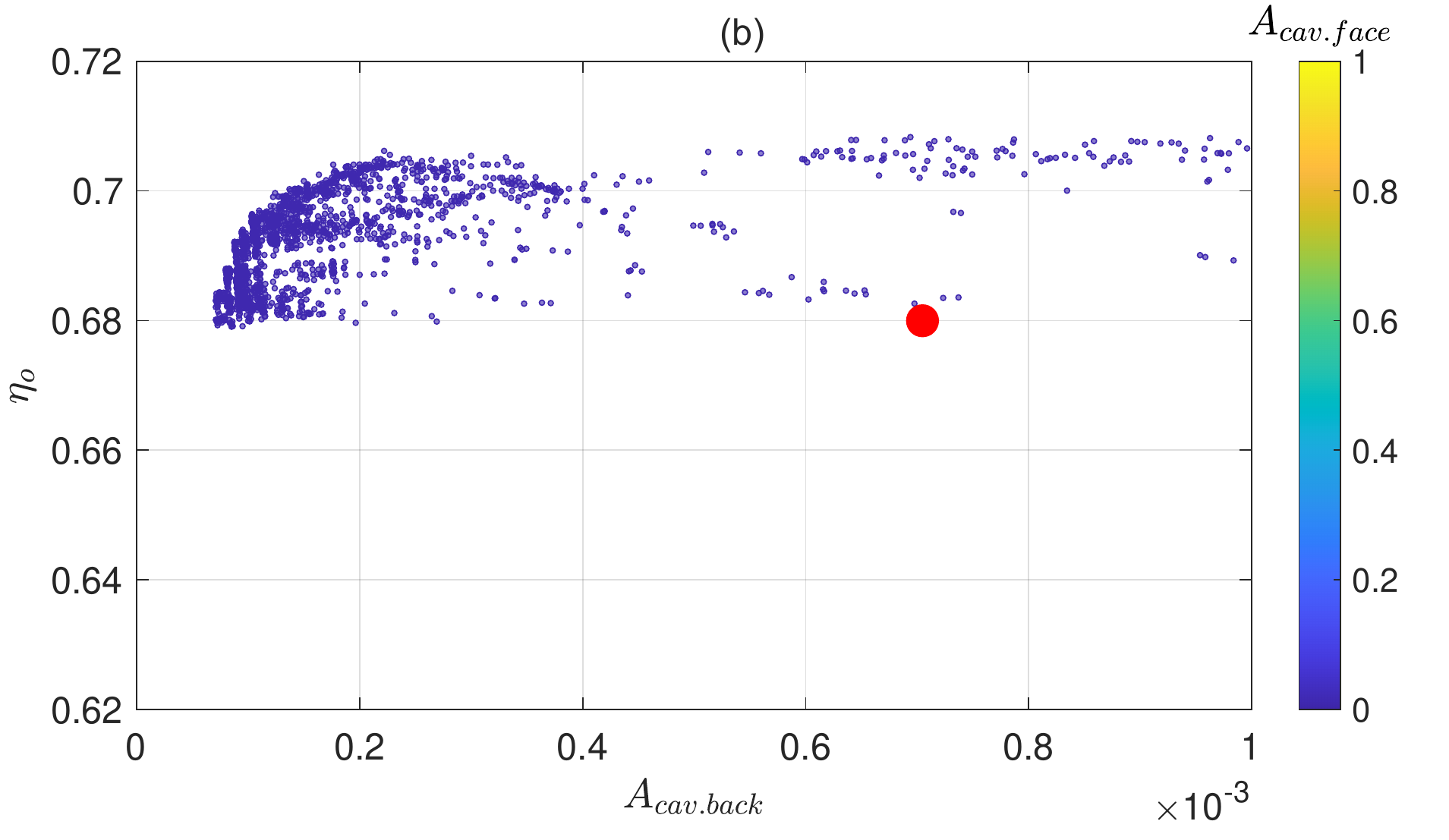}
    \caption{Pareto diagrams of the SDSO in the original non-reduced (a) and the SSDR (b) reduced spaces. The performances of the E779A selected test case are marked with a red dot.}
    \label{Feasible_NonReduced_s1}
    \end{figure}

    \subsubsection{Non-reduced design space}
We considered, at first, the simulation-driven shape optimisation problem in the non-reduced design space, i.e. using the 40-dimension space represented by the parameters of the B-Spline curves describing the primary radial and section geometrical characteristics of the blade from which all the samples for dimensions space reduction were derived. The results of this design activity represent the reference to assess the performances of the space dimensionality reduction applied to the design by optimization processes. 
To this aim, we filled the design space with 800 samples distributed accordingly to Uniform Latin-Hypercube bounded by the range of variation of the design parameters, which are similar to those adopted in \cite{gaggero2020reduced}. This initial population was allowed to evolve for 40 generations for a total of 32000 tested geometries. For each of them, all the Key Performance Indicators of the optimisation problem of Eq.\eqref{opt_eq}  were collected and organised in the Pareto diagram of Fig.\ref{Feasible_NonReduced_s1} (a). Cavitation was monitored separately on the suction ($A_{cav. back}$) and on the pressure ($A_{cav. face}$) sides of the blade since for the final selection of the optimal geometries we requested no cavitation on the face side.

As shown in several design activities concerning the E779A propeller \cite{gaggero2020reduced, GAGGERO_multifid}, there is the opportunity to completely nullify the sheet cavitation on the suction side ensuring, at the same time, no cavitation on the face of the blade. This result can be achieved at unvaried efficiency while, at the cost of a certain suction side cavitation (in any case reduced more than 85\% compared to the reference geometry), the efficiency can reach a 2.5\% increase.

    \subsubsection{Reduced design space}
    
\begin{figure}[htb!] 
    \centering
    \includegraphics[width=0.48\textwidth]{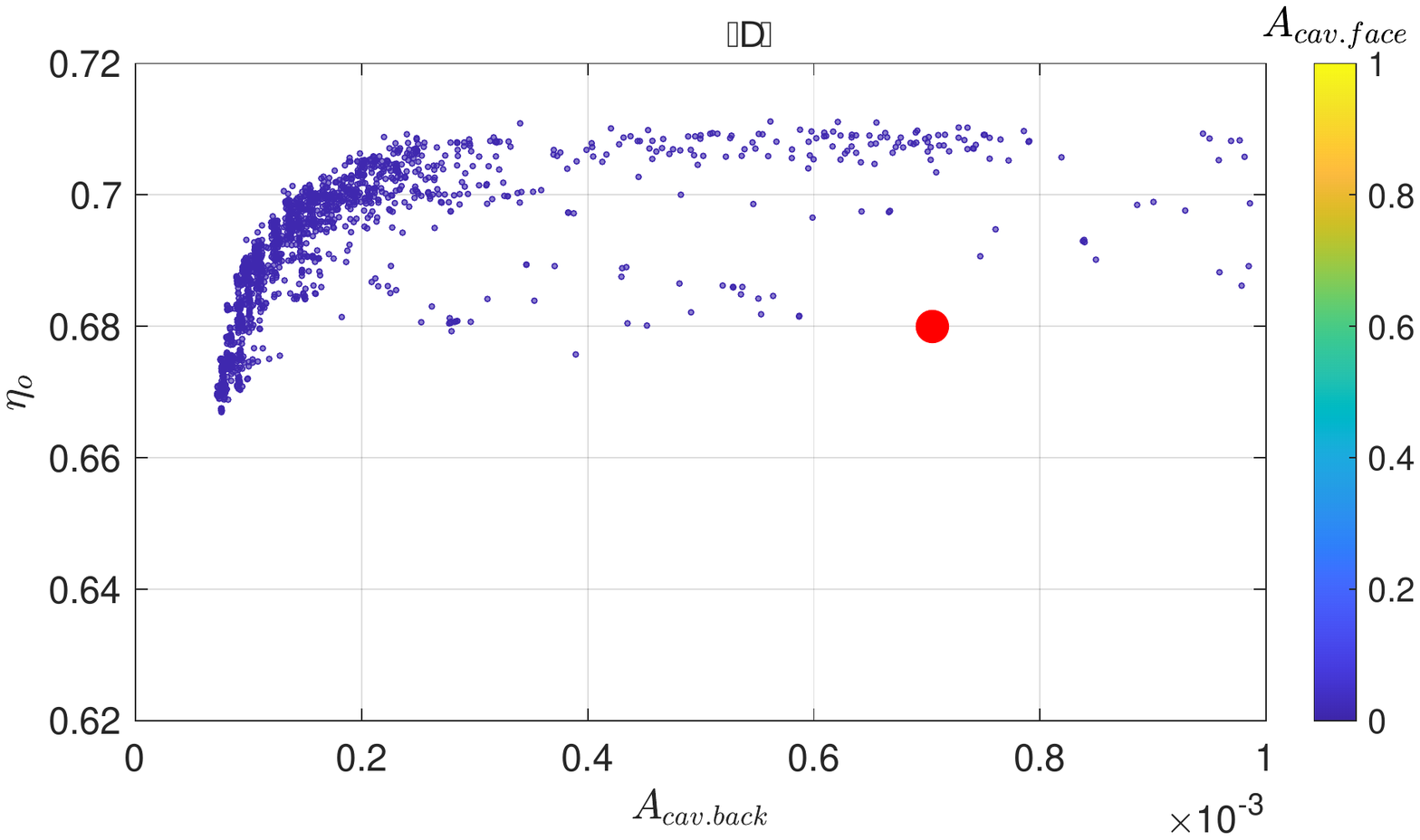}
    \includegraphics[width=0.48\textwidth]{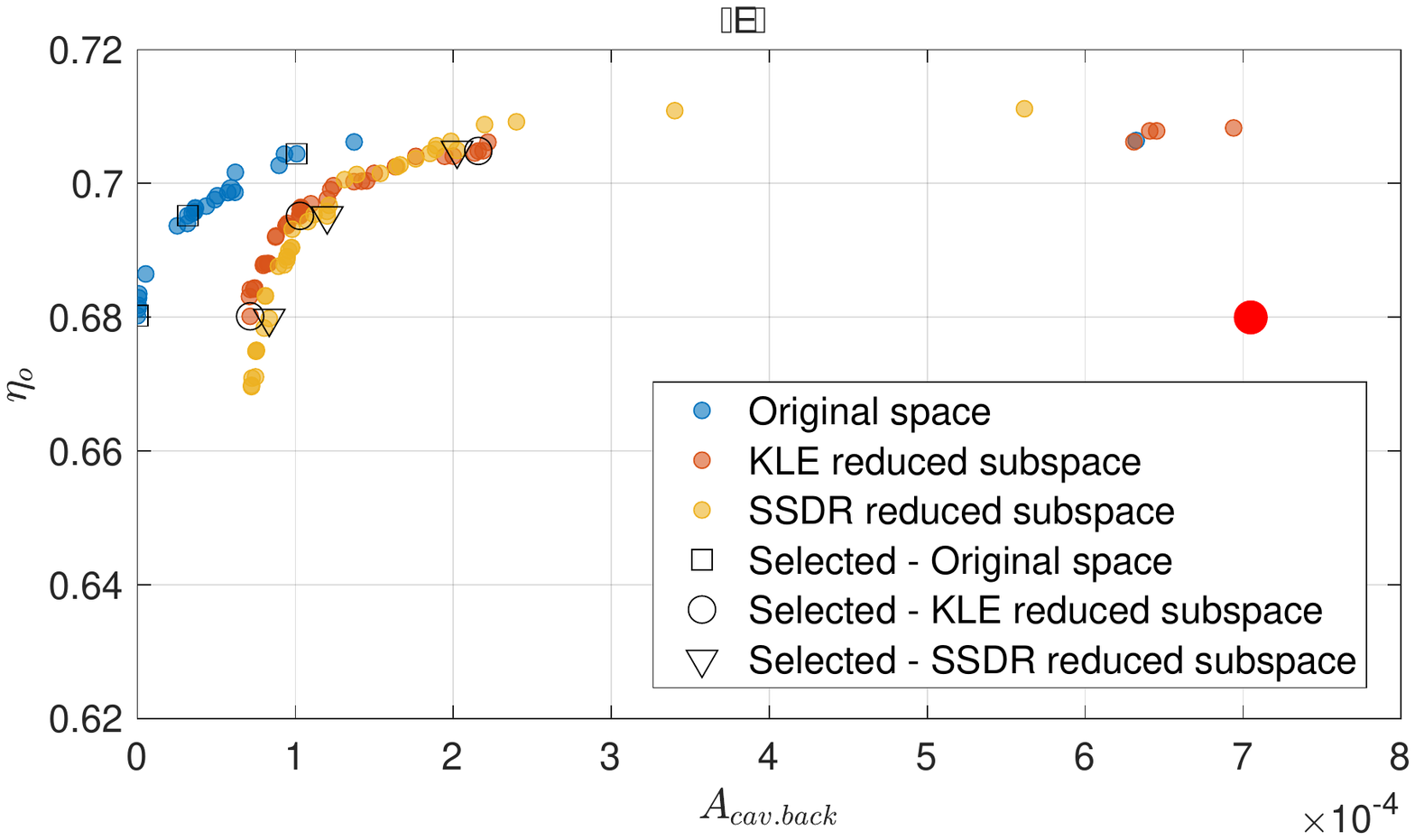}
    \caption{Pareto diagrams of the SDSO in the KLE reduced space (a) and comparison (b) of Pareto frontiers of the optimization activities. The performances of the E779A selected test case are marked with a red dot.}
    \label{Feasible_S1_cfr}
    \end{figure}

The simulation-driven shape optimisations using these reduced spaces require then an initial population of only 150 and 180 elements, respectively for SSDR and KLE ($10 \times m \times objectives$) distributed again using a Uniform Latin hypercube sampling approach. A total of 30 evolution of these initial populations was considered to collect converged solutions, evaluating a total of 4500 (SSDR) and 5400 (KLE) different geometries. Feasible designs, i.e., those satisfying the constraints on the delivered thrust, are collected in Fig. \ref{Feasible_NonReduced_s1} (b) (SSDR) and in Fig. \ref{Feasible_S1_cfr} (a) (KLE). 

SDSO in both spaces is successful. Compared to the performances of the reference E779A propeller, the design process using reduced spaces identifies geometries capable of providing a substantial reduction ($\approx 85\%$ at the efficiency of the reference propeller) of the cavitation extension on the back side (cavitation on the pressure side is always avoided) simultaneously with an increase of the propulsive efficiency that, similarly to the non-reduced space, reaches a 2\% of increment. 

    \begin{figure}[htb!] 
    \centering
    \includegraphics[width=0.3\textwidth]{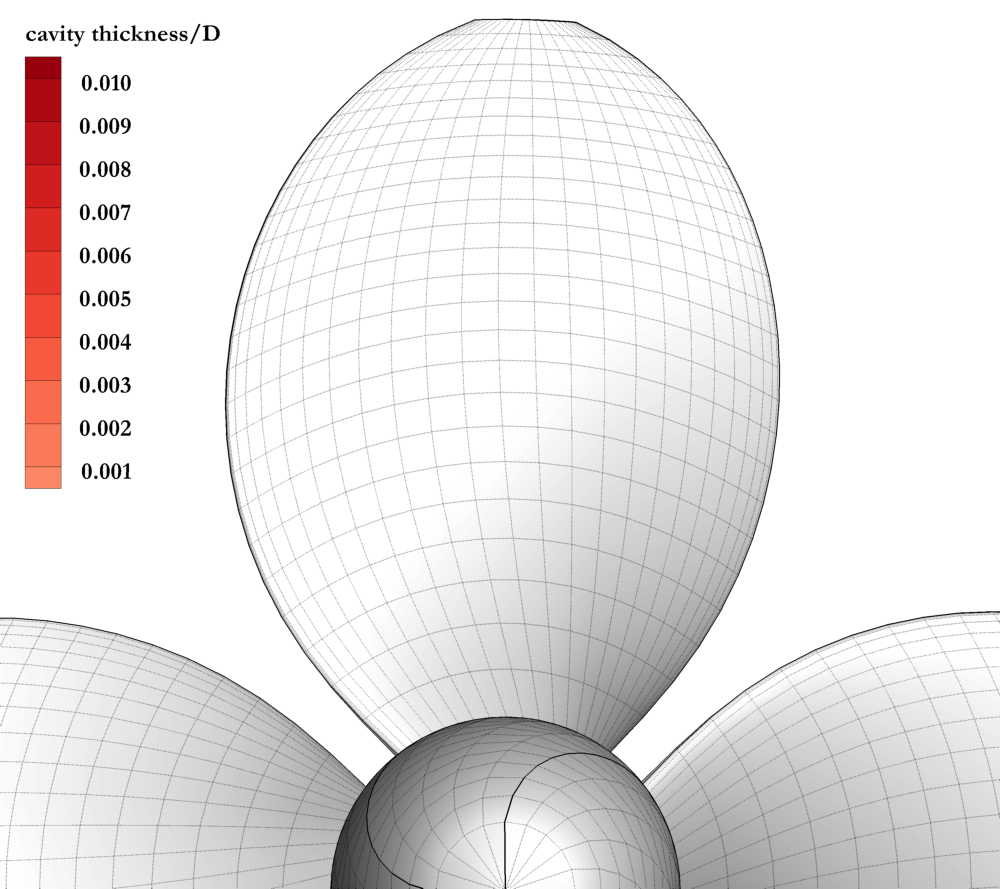}
    \includegraphics[width=0.3\textwidth]{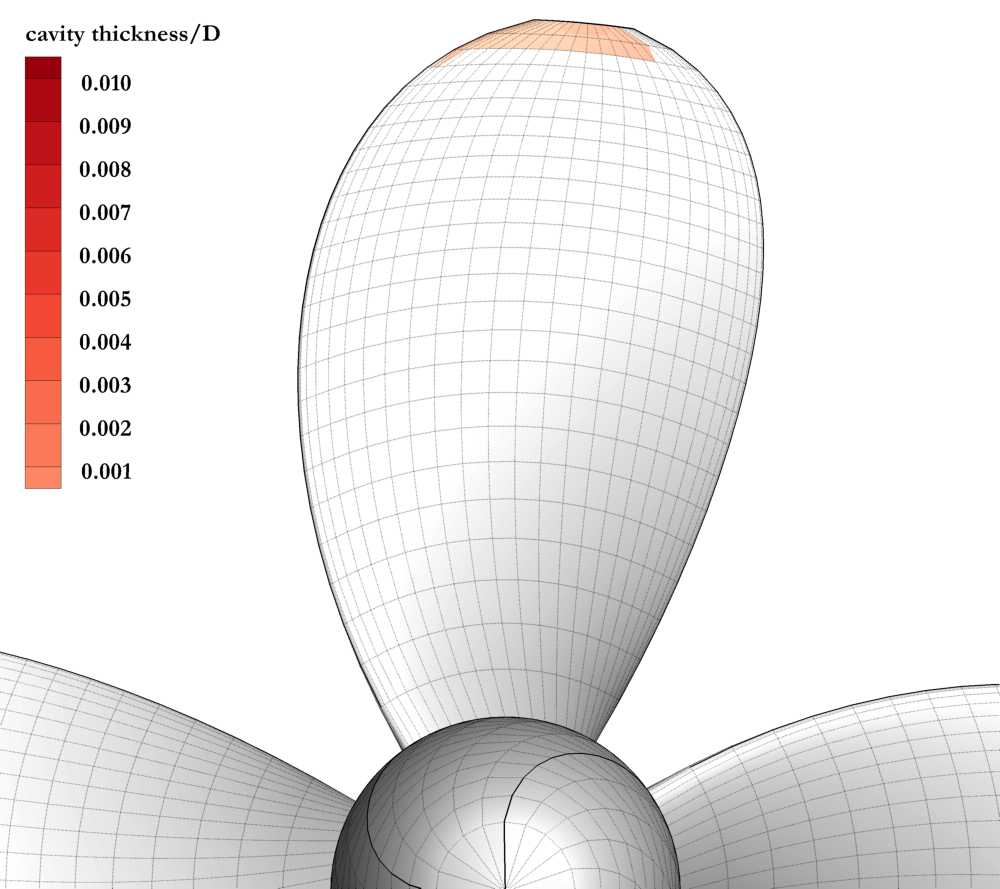}
    \includegraphics[width=0.3\textwidth]{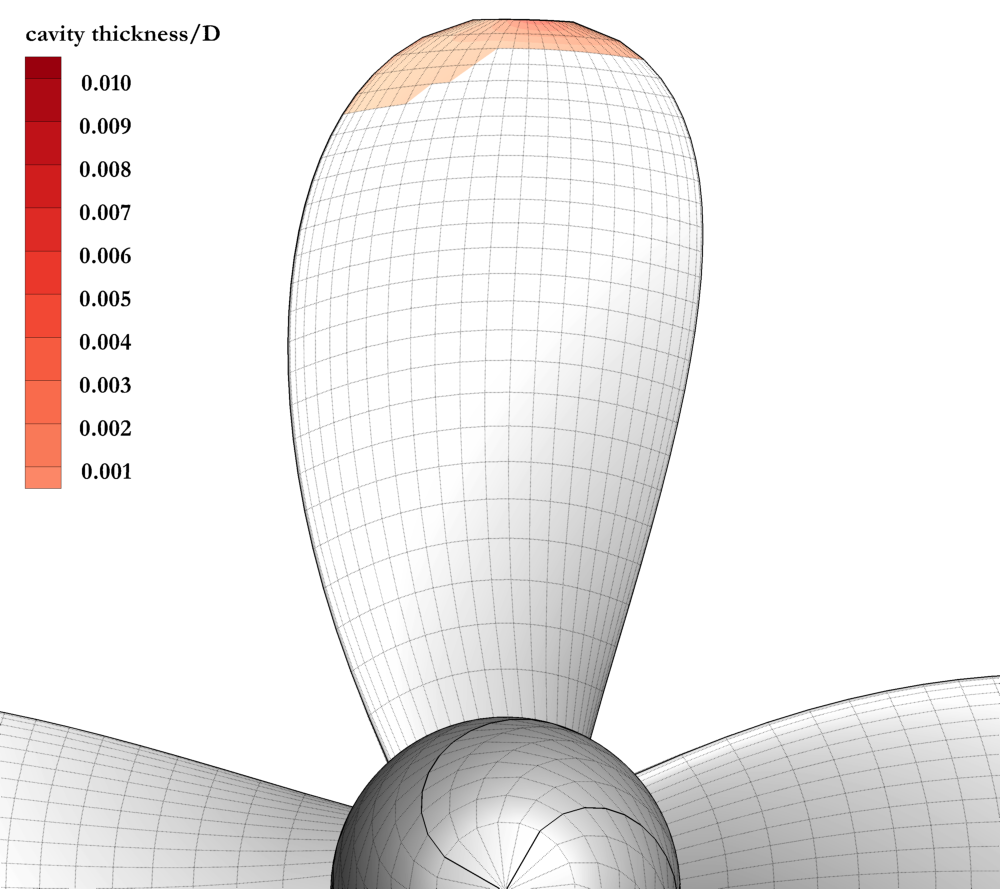}

    \includegraphics[width=0.3\textwidth]{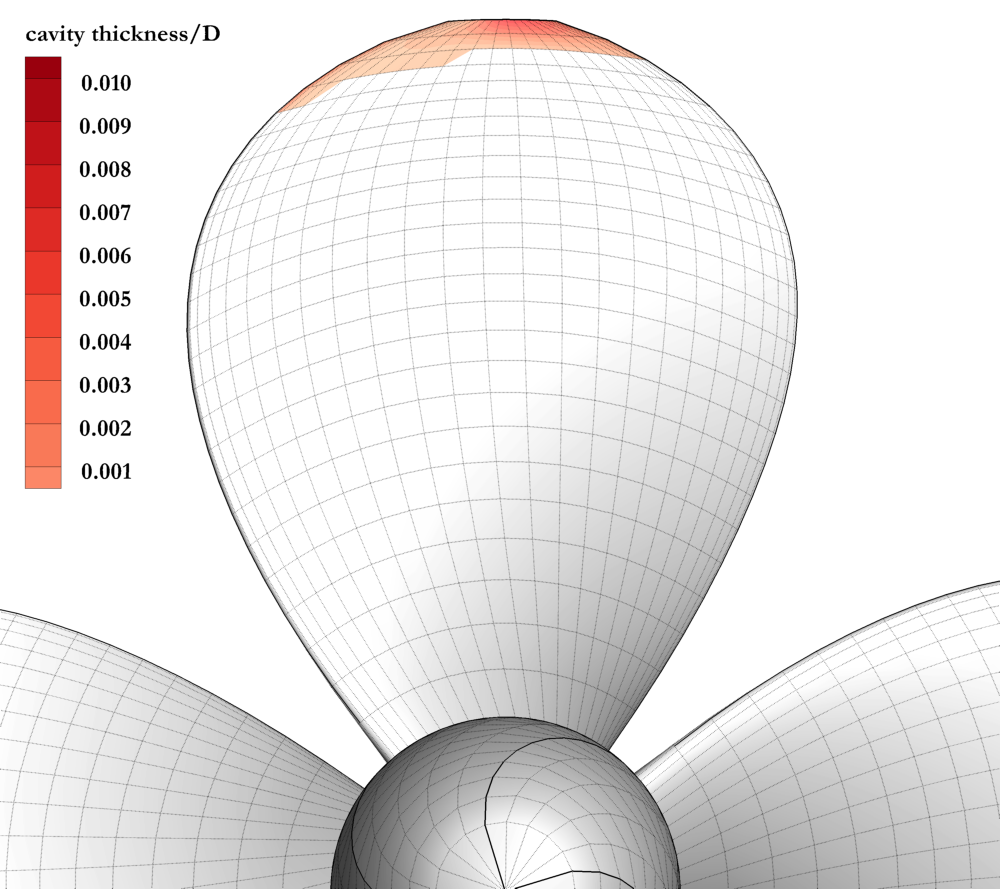}
    \includegraphics[width=0.3\textwidth]{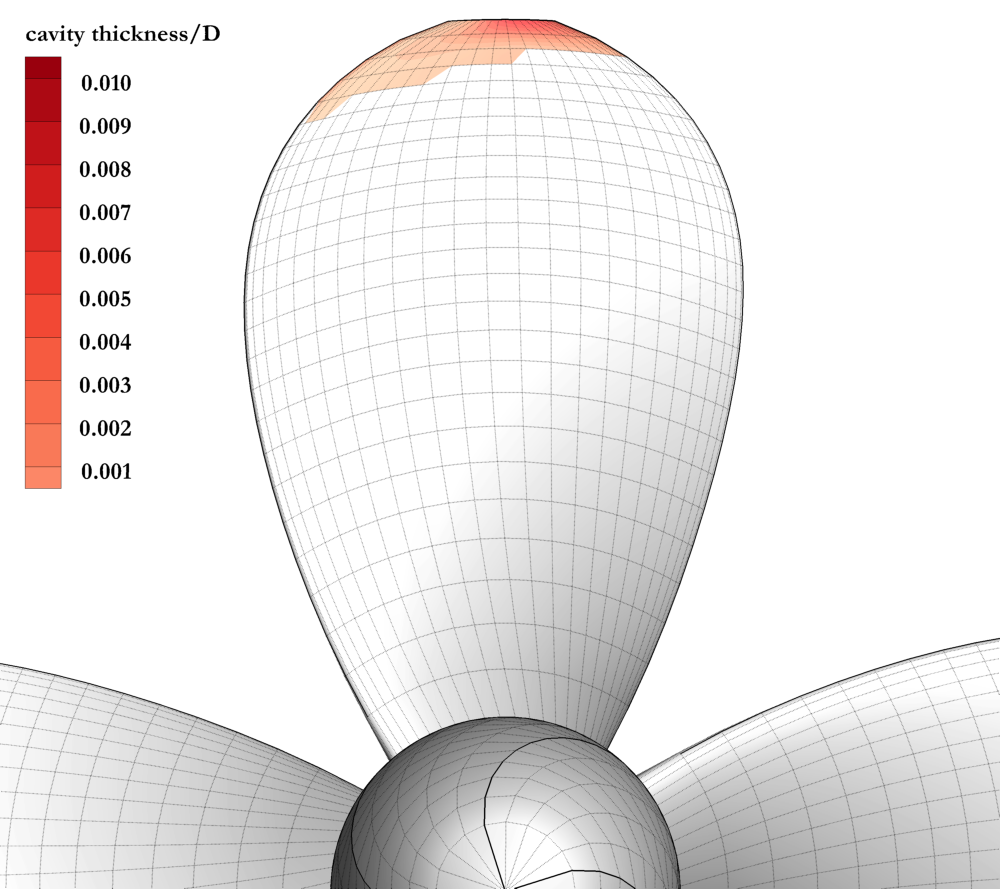}
    \includegraphics[width=0.3\textwidth]{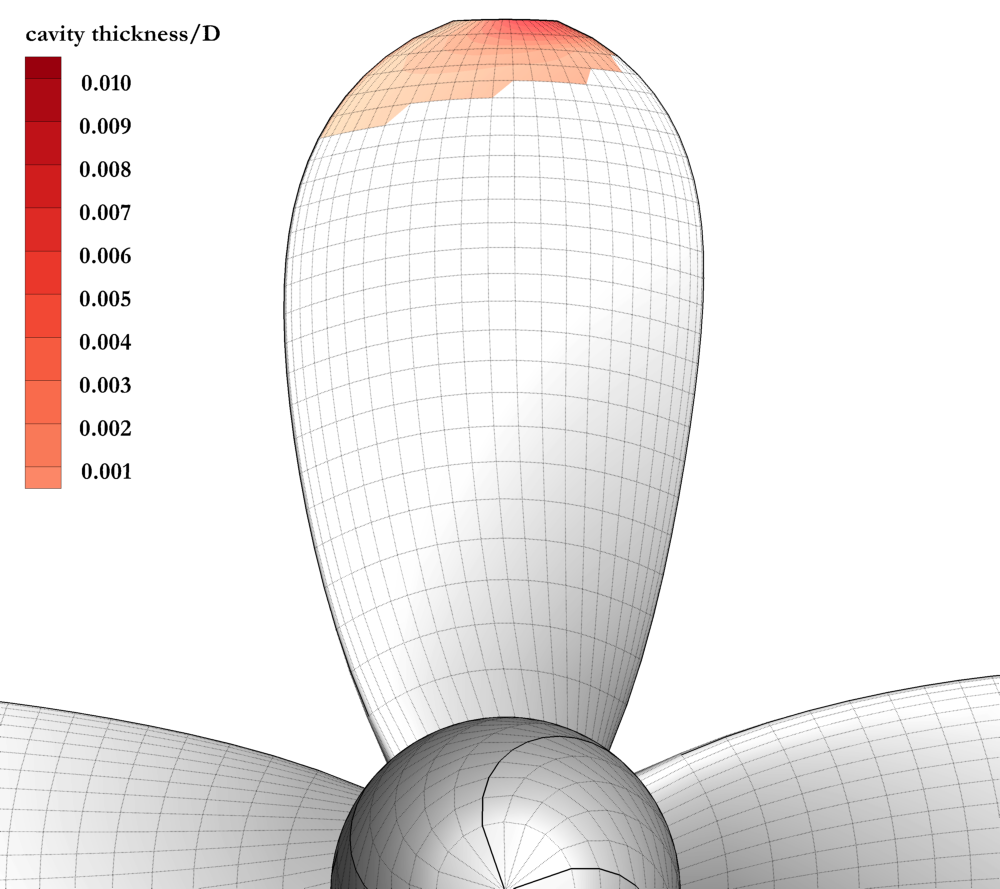}

    \includegraphics[width=0.3\textwidth]{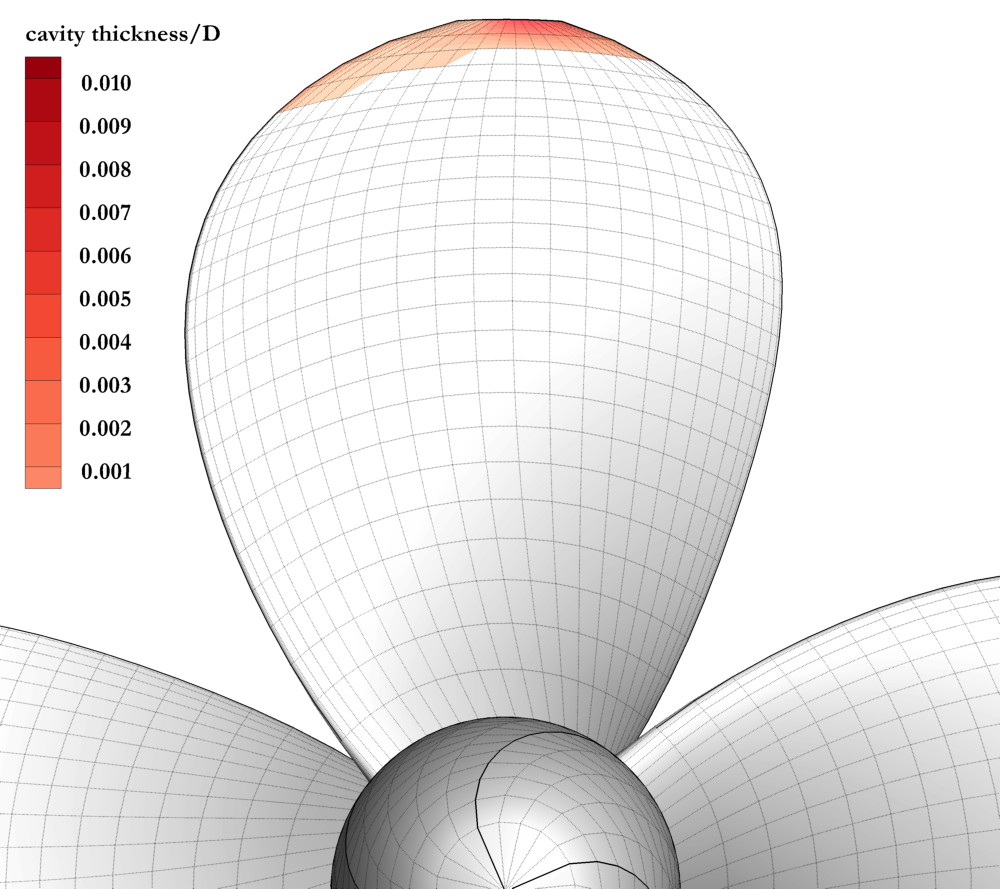}
    \includegraphics[width=0.3\textwidth]{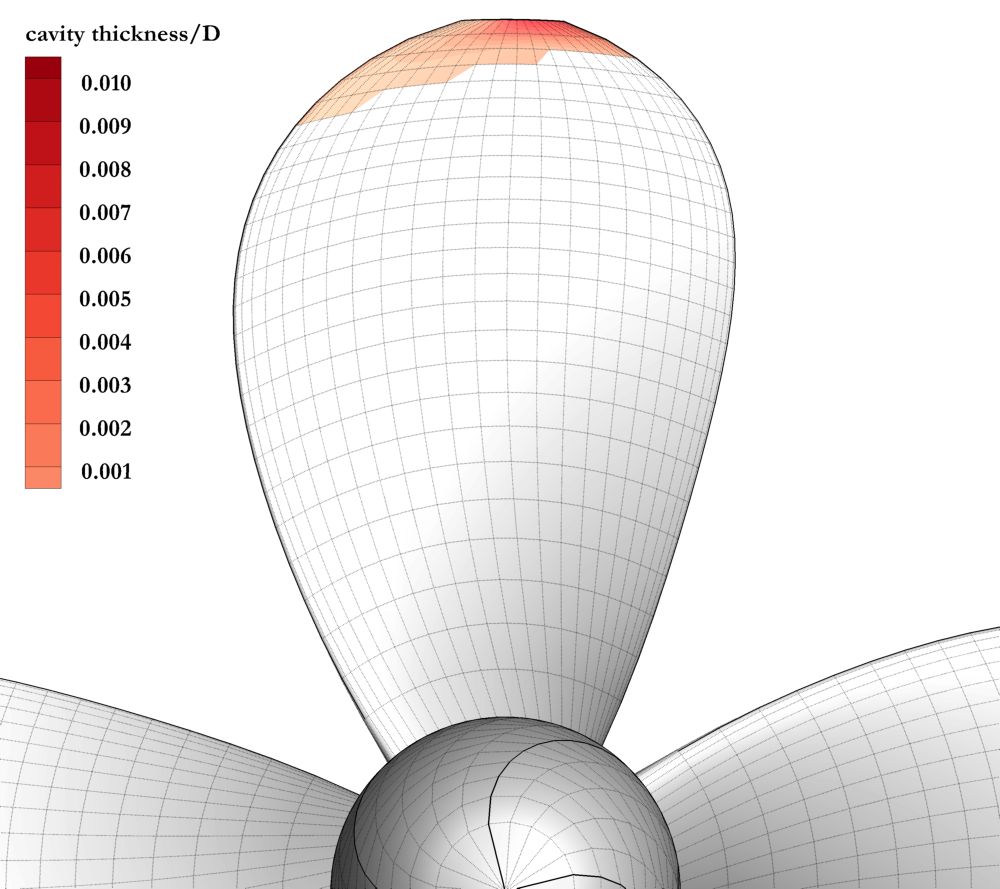}
    \includegraphics[width=0.3\textwidth]{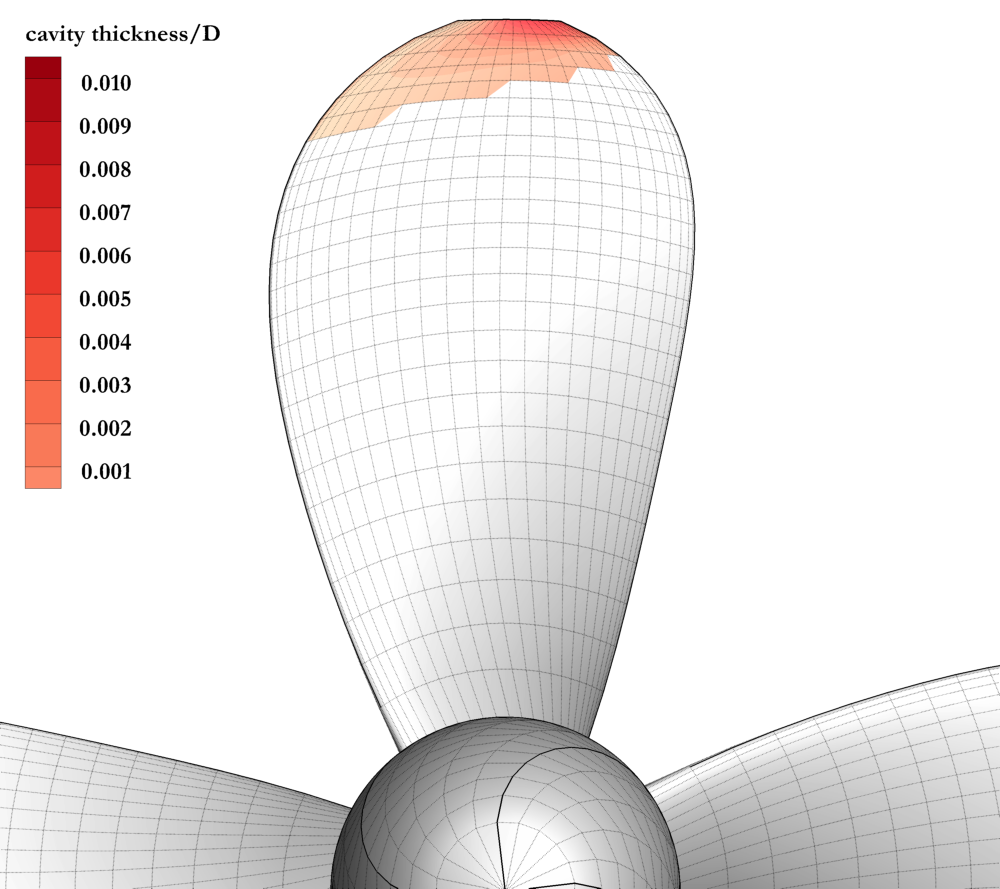}
    
    \caption{Performance comparison of selected Pareto geometries. From top to bottom: original, non-reduced space, KLE reduced space and SSDR reduced space. From left to right: efficiency level of 0.680, 0.695 and 0.705.}
    \label{Cfr_Cav}
    \end{figure}

Compared to the design process built on the non-reduced space, however, there are some differences that can be better appreciated considering the comparison of Pareto frontiers of Fig. \ref{Feasible_S1_cfr} (b). If in terms of optimisation, SSDR and KLE behave very similarly, and the addition of moments in the case of SSDR allows replicating the results of the larger KLE design space using fewer dimensions, from an absolute point of view the reduced spaces show certain limits, especially for what concern the avoidance of cavitation. None of the optimal geometries devised by both KLE and SSDR permits complete nullify the occurrence of back cavitation (even at very reduced values of efficiency). This suggests possible limitations of the current space dimensionality reduction process that, even neglecting only a few percentages of variance in the geometry, in the end, does not account for very local and small modifications of the blade shape (like those due to the camber) that instead are of fundamental importance for the avoidance of suction peaks and, then, for the occurrence of cavitation. The absence of pressure side cavitation on all the feasible geometries of the SDSO in the reduced spaces indirectly confirms this perception: none of the geometries under investigation has locally a sufficiently low pitch or an excessively high camber, which balanced combination would have resulted instead in perfectly shock-free blade geometries.

    \begin{figure}[htb!] 
    \centering
    \includegraphics[width=0.15\textwidth]{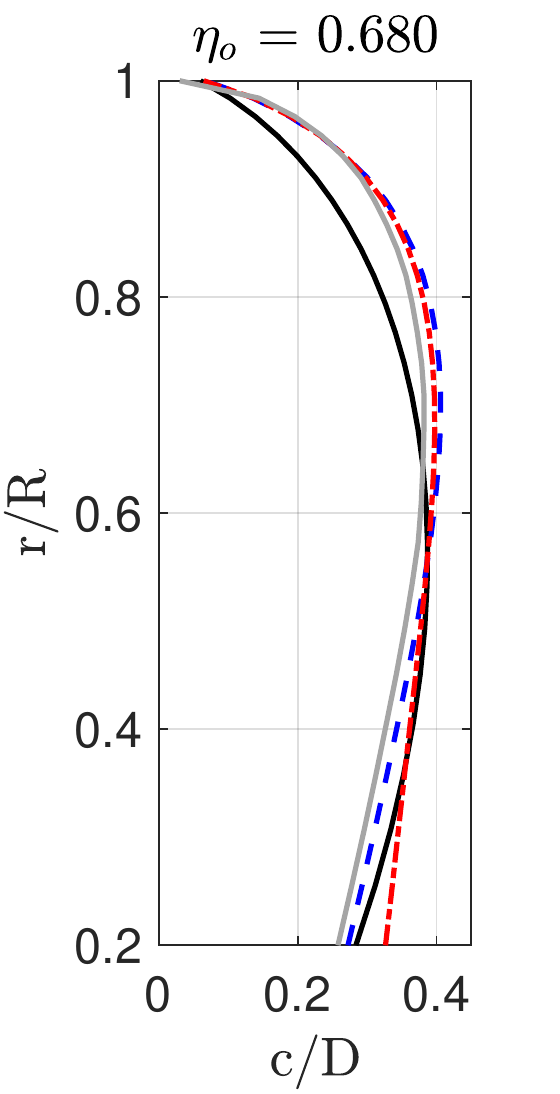}
    \includegraphics[width=0.15\textwidth]{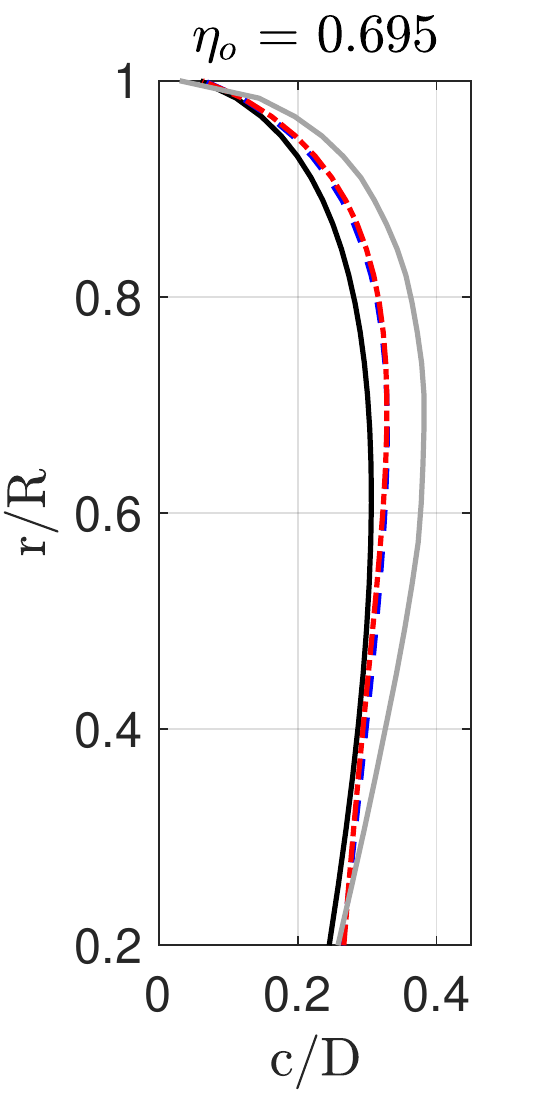}
    \includegraphics[width=0.15\textwidth]{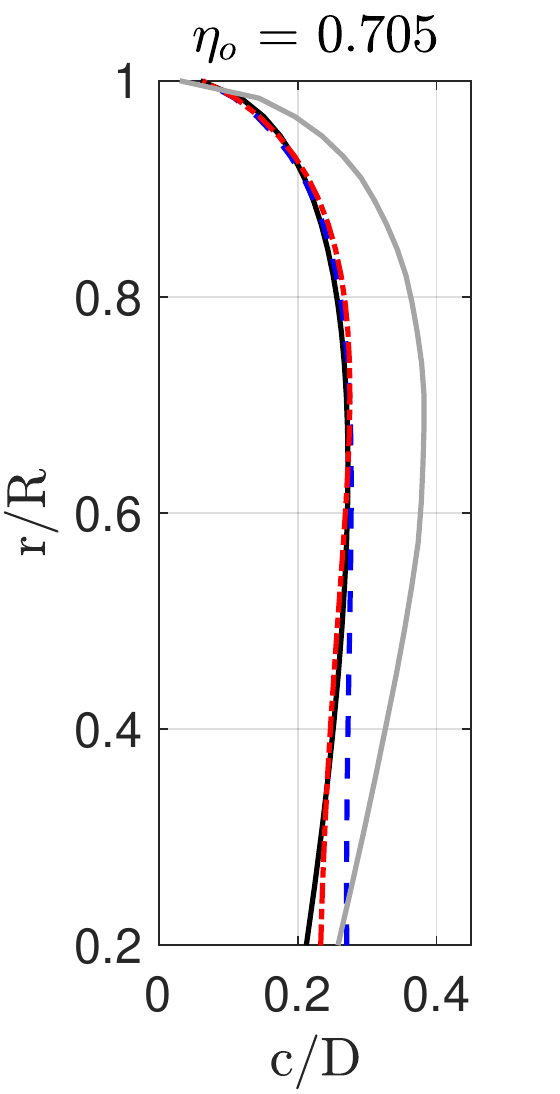}
    \includegraphics[width=0.14\textwidth]{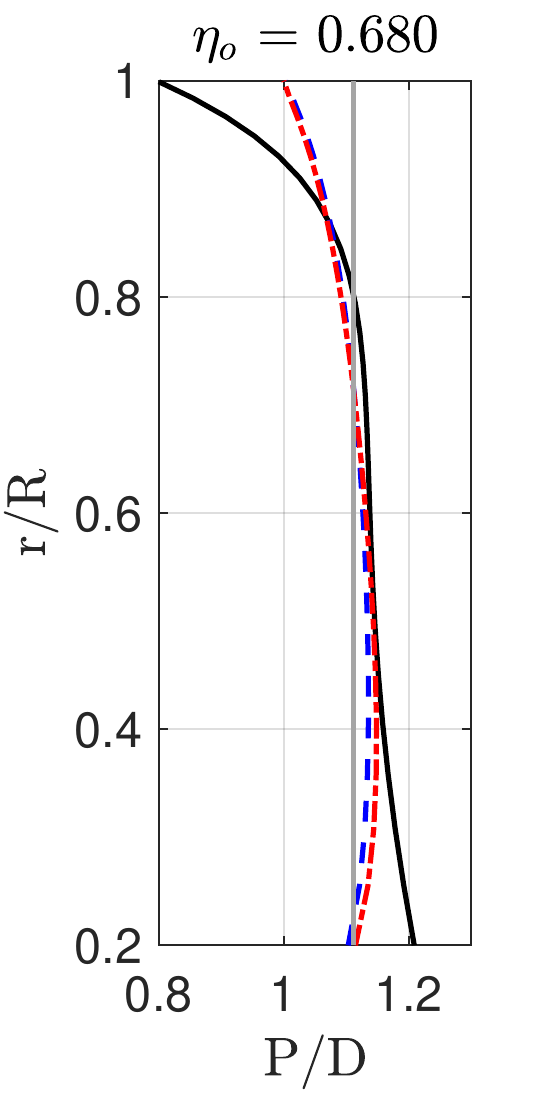}
    \includegraphics[width=0.14\textwidth]{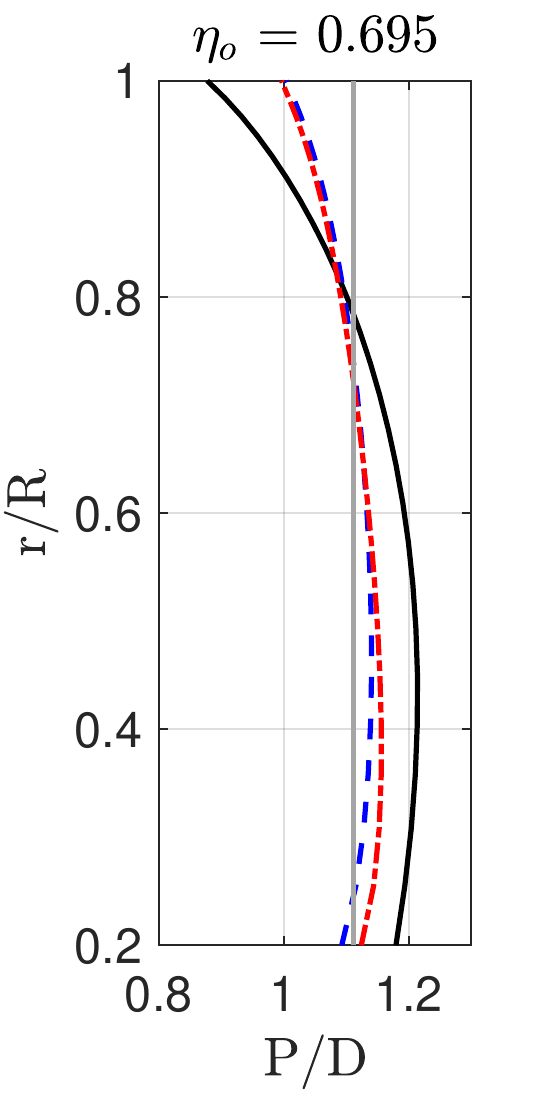}
    \includegraphics[width=0.14\textwidth]{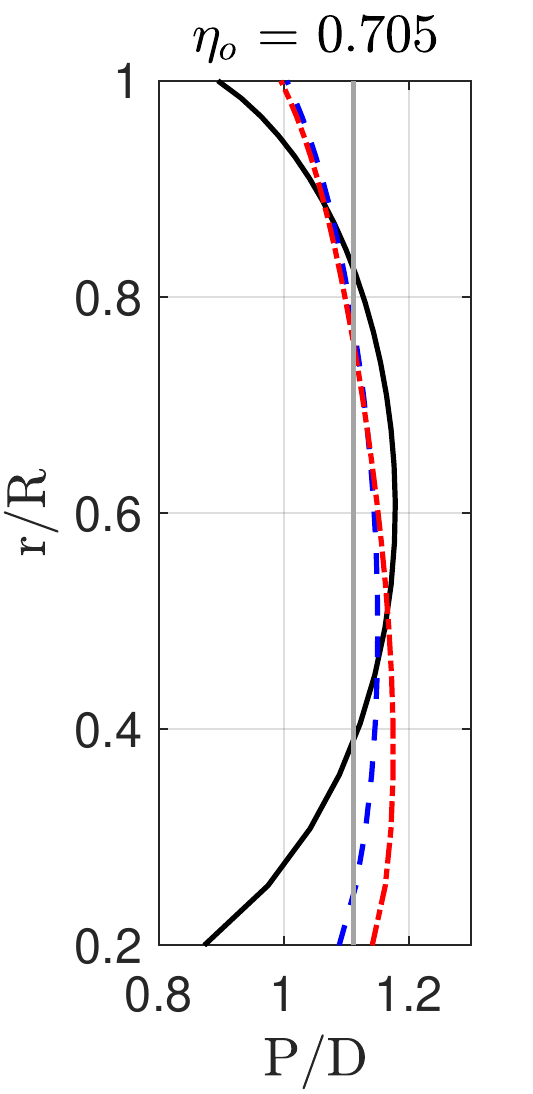}

    \includegraphics[width=0.17\textwidth]{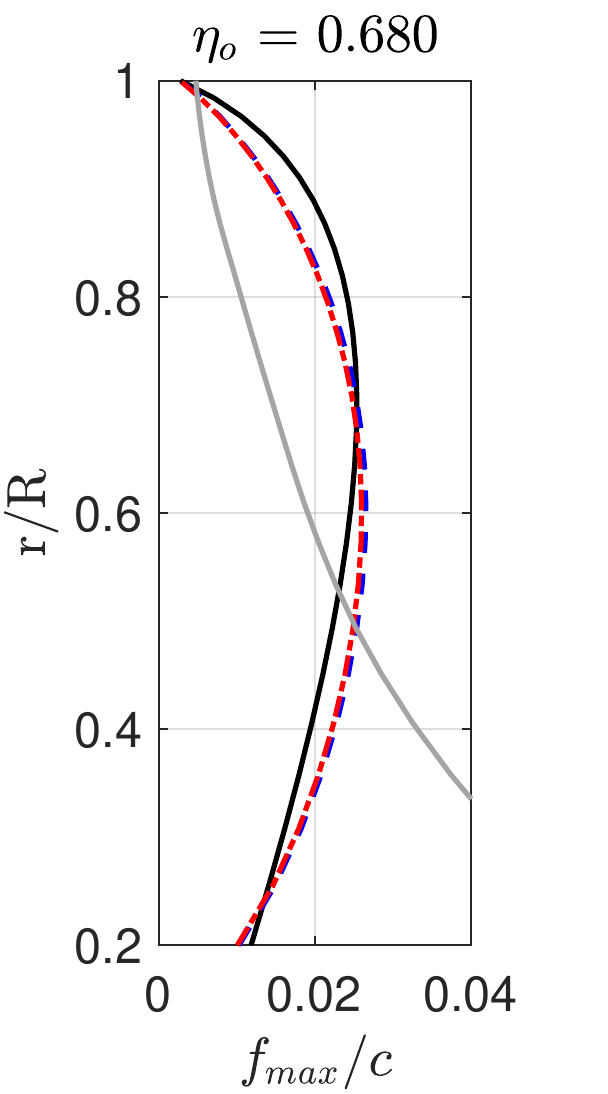}
    \includegraphics[width=0.17\textwidth]{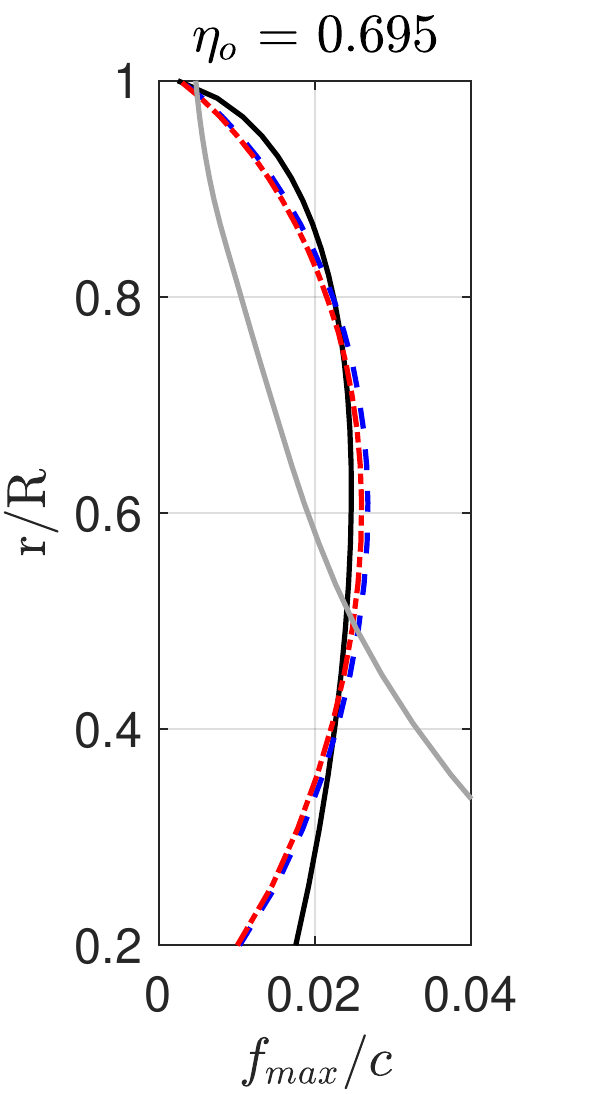}
    \includegraphics[width=0.17\textwidth]{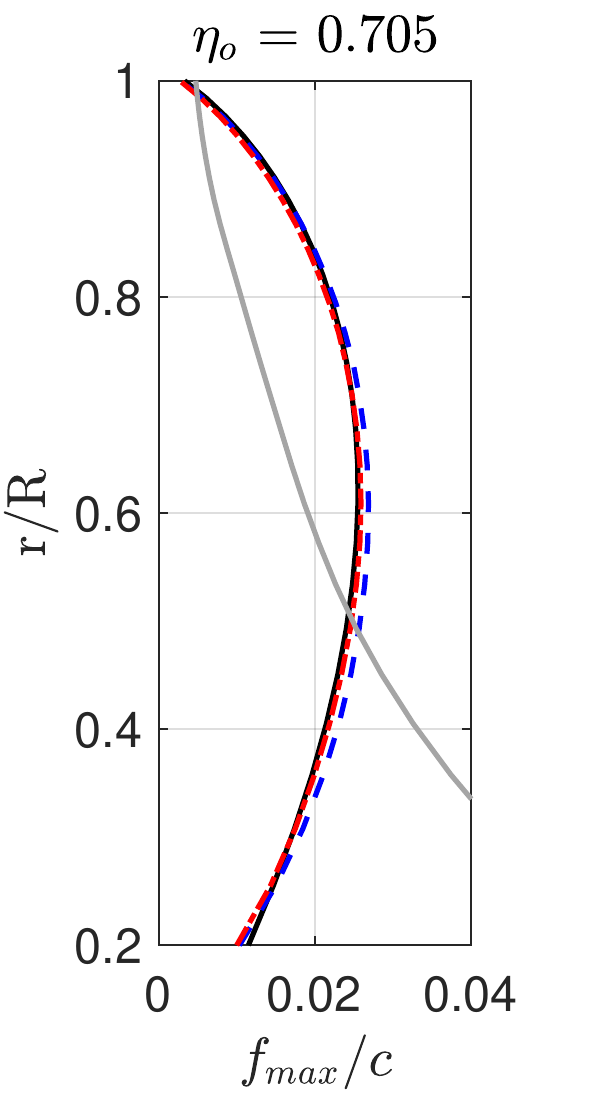}
    
    \caption{Comparison of selected Pareto geometries at different efficiency levels (0.680, 0.695 and 0.705). Non-dimensional chord, pitch and maximum sectional camber for non-reduced space (continuous black line), KLE (dashed blue) and SSDR (red, dash-dotted) reduced spaces. In grey is the geometry of the reference E779A propeller.}
    \label{Cfr_Geom}
    \end{figure}

To further discuss the results of the optimisation activities, based on the Pareto frontiers of Fig. \ref{Feasible_S1_cfr} (b), we selected a certain number of optimal geometries at constant propulsive efficiency (respectively, 0.680, 0.695 and 0.705 on the three design spaces considered) which performances in terms of cavitation extension on the back side, and of principal geometrical characteristics, are given in Fig. \ref{Cfr_Cav} and Fig. \ref{Cfr_Geom}. The effectiveness of both the reduced spaces for improving the cavitating performances of the propeller is clear. At the lowest efficiency, that of the reference propeller, as observed, the SDSO applied to the non-reduced design space identifies a completely cavitation-free geometry, but also the performance improvement of the geometries from KLE and SSDR is impressive if compared to the initial extension of cavitation on the reference propeller. At higher efficiencies, performances and propeller geometries among the spaces are closer and the cavitation provided by the optimal geometries from the reduced spaces is only slightly more extended than that from the non-reduced space at an efficiency level which is close to the maximum reachable in the entire design spaces. Geometry comparison highlights the reasons for these improvements. By looking at the chord distribution, we observe that, as expected, geometries providing higher values of efficiency have smaller expanded area ratio (i.e., shorter chords) and this trend is clear for geometries from the non-reduced design space as well as for geometries from KLE and SSDR. At the highest efficiency, chord distributions are almost overlapped, as also observable from the identical blade shapes of Fig. \ref{Cfr_Cav}.

Pitch and camber distributions, on the contrary, evidence the highest differences among the spaces. To reduce (and nullify) the extension of the cavity bubble, tip unloading is the easiest way. This is what naturally a designer would have done and this is what the SDSO applied on the non-reduced space pursues. The lowest cavitation (at the lowest efficiency) is achieved by a very reduced value of pitch at the tip that is never realised by the SDSO in the reduced spaces and this clarifies the presence of tip cavitation on the corresponding geometries. An increase of efficiency corresponds, in the original space, to an increase of pitch at the tip (i.e., tip loading) up to values closer to those of the reduced space geometries, explaining the similarities among Pareto configurations at higher efficiencies. We observe that the minimum value of pitch at the tip in the case of the KLE and SSDR is almost unchanged among spaces and geometries as if the minimum value possible in the spaces has been reached, and this represents a possible limitation of the space dimensionality reduction process. 
A similar discussion applies to the maximum sectional camber, which is very similarly realised regardless of the reduced spaces and the performances of the propeller (i.e., the efficiency level of the selection). Also in this case, reduced spaces geometries are very close to the non-reduced one at the highest efficiency, explaining the very similar performances observed in terms of cavitation extension.


\section{CONCLUSIONS}

In this work, we use the shape-supervised DR technique, which uses higher-level information for the shape in terms of its geometric integral properties. In our case, these integral properties are geometric moments of varying order, whose usage is based on two fundamental insights: (i) the geometric moments of a shape are intrinsic features of the underlying geometry and (ii) they provide a unifying medium between the shape's geometry and physics. To maximise the geometric and physical information retained in the subspace, we evaluate geometric moments using the divergence theorem. Once moments are evaluated, they are used, along with the shape modification function, to form a shape signature vector (SSV), which acts as a descriptor to uniquely represent each instance in the design space. Afterwards, we use SSV to construct a symmetric and positive definite covariance matrix, whose eigendecomposition results in a latent feature matrix. The columns of this matrix are orthogonal eigenvectors, which, along with the highest eigenvalues, span the basis of a subspace retaining the highest variance in terms of geometry, its underlying structure and physics. 

We used the E779A propeller for the experimentation in this work, which was parameterised with 40 design parameters defining the primary geometrical characteristics of the propeller, including pitch, chord, maximum and sectional camber distributions along the blade radius. The shape of the propeller was optimised to maximise efficiency while reducing suction side cavitation.

The subspace achieved using both KLE and SSDR approaches results in similar levels of dimensionality reduction. However, the subspace obtained from the SSDR approach is of better quality than that obtained from the KLE approach, as it provides a significantly higher number of valid geometries. The optimisation is then performed first in the original 40-dimensional design space and subsequently in the subspaces generated by the KLE and SSDR approaches. Due to the high dimensionality, convergence in the original design space is significantly slower compared to the subspaces. However, there is no significant difference in performance between designs obtained from the original design space and the reduced-dimensional subspaces.

\section*{ACKNOWLEDGEMENTS}
This work received funding from the European Union's Horizon 2020 research and innovation programme under the Marie Skłodowska-Curie grant agreement No 860843, PI for the University of Strathclyde: P.D. Kaklis.

\bibliographystyle{elsarticle-num}
\bibliography{ref}

\begin{thebibliography}{10}
\expandafter\ifx\csname url\endcsname\relax
  \def\url#1{\texttt{#1}}\fi
\expandafter\ifx\csname urlprefix\endcsname\relax\def\urlprefix{URL }\fi
\expandafter\ifx\csname href\endcsname\relax
  \def\href#1#2{#2} \def\path#1{#1}\fi

\bibitem{gaggero2020reduced}
S.~Gaggero, G.~Vernengo, D.~Villa, L.~Bonfiglio, A reduced order approach for
  optimal design of efficient marine propellers, Ships and Offshore Structures
  15~(2) (2020) 200--214.

\bibitem{khan2017customer}
S.~Khan, E.~Gunpinar, M.~Moriguchi, Customer-centered design sampling for {CAD}
  products using spatial simulated annealing, Proceedings of CAD 17 (2017)
  100--103.

\bibitem{R07}
M.~Diez, E.~F. Campana, F.~Stern, Design-space dimensionality reduction in
  shape optimization by karhunen--lo{\`e}ve expansion, Computer Methods in
  Applied Mechanics and Engineering 283 (2015) 1525--1544.

\bibitem{d2017nonlinear}
D.~D’Agostino, A.~Serani, E.~F. Campana, M.~Diez, Nonlinear methods for
  design-space dimensionality reduction in shape optimization, in:
  International Workshop on Machine Learning, Optimization, and Big Data,
  Springer, 2017, pp. 121--132.

\bibitem{masood2021machine}
Z.~Masood, S.~Khan, L.~Qian, Machine learning-based surrogate model for
  accelerating simulation-driven optimisation of hydropower kaplan turbine,
  Renewable Energy 173 (2021) 827--848.

\bibitem{khan2022shape}
S.~Khan, P.~Kaklis, A.~Serani, M.~Diez, K.~Kostas, Shape-supervised dimension
  reduction: Extracting geometry and physics associated features with geometric
  moments, Computer-Aided Design 150 (2022) 103327.

\bibitem{KHAN2023116051}
S.~Khan, K.~Goucher-Lambert, K.~Kostas, P.~Kaklis, Ship{H}ull{GAN}: A generic
  parametric modeller for ship hull design using deep convolutional generative
  model, Computer Methods in Applied Mechanics and Engineering 411 (2023)
  116051.

\bibitem{khan2022geometric}
S.~Khan, P.~Kaklis, A.~Serani, M.~Diez, Geometric moment-dependent global
  sensitivity analysis without simulation data: application to ship hull form
  optimisation, Computer-Aided Design 151 (2022) 103339.

\bibitem{khan2023shiphullgan}
S.~Khan, K.~Goucher-Lambert, K.~Kostas, P.~Kaklis, Shiphullgan: A generic
  parametric modeller for ship hull design using deep convolutional generative
  model, Computer Methods in Applied Mechanics and Engineering 411 (2023)
  116051.

\bibitem{salvatore2006description}
F.~Salvatore, F.~Pereira, M.~Felli, D.~Calcagni, F.~Di~Felice, Description of
  the insean e779a propeller experimental dataset, Technical Report INSEAN
  2006-085 INSEAN-Italian Ship Model Basin (2006).

\bibitem{bertetta2012cpp}
D.~Bertetta, S.~Brizzolara, S.~Gaggero, M.~Viviani, L.~Savio, Cpp propeller
  cavitation and noise optimization at different pitches with panel code and
  validation by cavitation tunnel measurements, Ocean engineering 53 (2012)
  177--195.

\bibitem{gaggero2016design}
S.~Gaggero, J.~Gonzalez-Adalid, M.~P. Sobrino, Design of contracted and tip
  loaded propellers by using boundary element methods and optimization
  algorithms, Applied Ocean Research 55 (2016) 102--129.

\bibitem{pereira2002experimental}
F.~Pereira, F.~Salvatore, F.~Di~Felice, M.~Elefante, Experimental and numerical
  investigation of the cavitation pattern on a marine propeller, in: 24th
  Symposium on Naval Hydrodynamics. Fukuoka, Japan, 2002, pp. 8--13.

\bibitem{salvatore2009propeller}
F.~Salvatore, H.~Streckwall, T.~Van~Terwisga, Propeller cavitation modelling by
  cfd-results from the virtue 2008 rome workshop, in: Proceedings of the first
  international symposium on marine propulsors, Trondheim, Norway, 2009, pp.
  22--24.

\bibitem{bensow2010simulating}
R.~E. Bensow, G.~Bark, Simulating cavitating flows with les in openfoam, in: V
  European conference on computational fluid dynamics, 2010, pp. 14--17.

\bibitem{vaz2015cavitating}
G.~Vaz, D.~Hally, T.~Huuva, N.~Bulten, P.~Muller, P.~Becchi, J.~Herrer,
  S.~Whitworth, R.~Mac{\'e}, A.~Korsstr{\"o}m, Cavitating flow calculations for
  the e779a propeller in open water and behind conditions: code comparison and
  solution validation, in: Proceedings of the Fourth International Symposium on
  Marine Propulsors (SMP), Vol.~15, 2015, pp. 344--360.

\bibitem{yang1997fast}
L.~Yang, F.~Albregtsen, T.~Taxt, Fast computation of three-dimensional
  geometric moments using a discrete divergence theorem and a generalization to
  higher dimensions, Graphical models and image processing 59~(2) (1997)
  97--108.

\bibitem{M1}
A.~Krishnamurthy, S.~McMains, Accurate gpu-accelerated surface integrals for
  moment computation, Computer-Aided Design 43~(10) (2011) 1284--1295.

\bibitem{M2}
A.~Taber, G.~Kumar, M.~Freytag, V.~Shapiro, A moment-vector approach to
  interoperable analysis, Computer-Aided Design 102 (2018) 139--147.

\bibitem{M10}
D.~F. Atrevi, D.~Vivet, F.~Duculty, B.~Emile, A very simple framework for 3d
  human poses estimation using a single 2d image: Comparison of geometric
  moments descriptors, Pattern Recognition 71 (2017) 389--401.

\bibitem{M7}
A.~M. Bronstein, M.~M. Bronstein, R.~Kimmel, Numerical geometry of non-rigid
  shapes, Springer Science \& Business Media, 2008.

\bibitem{M8}
P.~Jin, B.~Xie, F.~Xiao, Multi-moment finite volume method for incompressible
  flows on unstructured moving grids and its application to fluid-rigid body
  interactions, Computers \& Structures 221 (2019) 91--110.

\bibitem{xu2008geometric}
D.~Xu, H.~Li, Geometric moment invariants, Pattern recognition 41~(1) (2008)
  240--249.

\bibitem{GAGGERO_multifid}
S.~Gaggero, G.~Vernengo, D.~Villa, A marine propeller design method based on
  two-fidelity data levels, Applied Ocean Research 123 (2022) 103156.

\end{thebibliography}
\end{document}